\documentclass[11pt]{article}

\usepackage{color}
\usepackage{latexsym}
\usepackage{dsfont}
\usepackage{amssymb}
\usepackage{graphicx}
\usepackage{amsmath, amsfonts,amssymb,theorem,euscript,array,enumerate,amsfonts,mathrsfs}
\usepackage{hyperref}

\newtheorem{Theorem}{Theorem}[part]
\newtheorem{Definition}{Definition}[part]
\newtheorem{Proposition}{Proposition}[part]

\newtheorem{Remark}{Remark}[part]

\def \trans{^{\scriptscriptstyle{\intercal}}}

\def \Frac{\displaystyle\frac}
\def \Inf{\displaystyle\inf}

\def \b1{\bf{1}}

\def \R{\mathbb{R}}

\def \E{\mathbb{E}}
\def \F{\mathbb{F}}

\def \H{\mathbb{H}}
\def \P{\mathbb{P}}

\def \S{\mathbb{S}}

\def\esssup_#1{\underset{#1}{\mathrm{ess\,sup\, }}}
\def\essinf_#1{\underset{#1}{\mathrm{ess\,inf\, }}}

\def\argmin_#1{\underset{#1}{\mathrm{argmin\, }}}

\def \Ac{{\cal A}}
\def \Bc{{\cal B}}
\def \Cc{{\cal C}}

\def \Fc{{\cal F}}

\def \Hc{{\cal H}}

\def \Kc{{\cal K}}
\def \Lc{{\cal L}}
\def \Pc{{\cal P}}

\def \Vc{{\cal V}}
\def \Wc{{\cal W}}

\def \eps{\varepsilon}

\def \ep{\hbox{ }\hfill$\Box$}

\def\Dt#1{\Frac{\partial #1}{\partial t}}

\def\reff#1{{\rm(\ref{#1})}}

\def\beqs{\begin{eqnarray*}}
\def\enqs{\end{eqnarray*}}
\def\beq{\begin{eqnarray}}
\def\enq{\end{eqnarray}}

\begin{document}

 \title{Bellman equation and viscosity solutions for \\ mean-field stochastic control problem
 %\thanks{We would like to thank E. Bayraktar, J.F. Chasagneux and A. Cosso for discussions during the preparation of this paper.}
 %\thanks{This study was supported by FiME (Finance for Energy Market Research Centre) and the ``Finance et D\'eveloppement Durable - Approches Quantitatives'' EDF - CACIB Chair.  }
 }

\author{Huy\^en PHAM\thanks{Corresponding author. Research  supported in part  by FiME (Finance for Energy Market Research Centre) and
the ``Finance et D\'eveloppement Durable - Approches Quantitatives'' EDF - CACIB Chair. }
\\\small  Laboratoire de Probabilit\'es et
 \\\small  Mod\`eles Al\'eatoires, CNRS, UMR 7599
 \\\small  Universit\'e Paris Diderot
 \\\small  Case courrier 7012
 \\\small  Avenue de France
 \\\small 75205 Paris Cedex 13, France
 \\\small  pham at math.univ-paris-diderot.fr
\\\small  and CREST-ENSAE
\and
Xiaoli WEI
\\\small  Laboratoire de Probabilit\'es et
 \\\small  Mod\`eles Al\'eatoires, CNRS, UMR 7599
 \\\small  Universit\'e Paris Diderot
  \\\small  Case courrier 7012
 \\\small  Avenue de France
 \\\small 75205 Paris Cedex 13, France
 \\\small  tyswxl at gmail.com
}

%\author{Huy{\^e}n PHAM\thanks{Laboratoire de Probabilit\'es et Mod\`eles Al\'eatoires, CNRS, UMR 7599, Universit{\'e} Paris Diderot, and
 %              CREST-ENSAE,  \sf pham at math.univ-paris-diderot.fr} ~~~
%		Xiaoli WEI\thanks{Laboratoire de Probabilit\'es et Mod\`eles Al\'eatoires, CNRS, UMR 7599, Universit\'e Paris Diderot,
%		\sf tyswxl at gmail.com}
 %            }

\maketitle

\date{}

\begin{abstract}
We consider the stochastic optimal control problem of  McKean-Vlasov stochastic differential equation where the coefficients may depend upon the joint law of the 
state and control.   By using feedback controls, we reformulate the problem into a deterministic control problem with only  the marginal distribution of the process as controlled state variable, 
and prove that dynamic programming principle holds in its general form.  
Then, by relying on the notion of differentiability with respect to pro\-bability measures recently introduced by P.L. Lions in \cite{lio12}, and a special 
It\^o formula for flows of  probability measures,  we derive the (dynamic programming) Bellman equation for mean-field stochastic control problem, and prove a veri\-fication theorem in our McKean-Vlasov framework.  We give explicit solutions to the Bellman equation for the linear quadratic mean-field control problem, with applications to the mean-variance portfolio selection and a systemic risk model.  We also consider a notion of lifted 
visco\-sity solutions for the Bellman equation, and show the viscosity pro\-perty and uniqueness of the value function to the McKean-Vlasov 
control problem. Finally, we consider the case of McKean-Vlasov control problem with open-loop controls and discuss the asso\-ciated dynamic programming equation that we compare with the case of closed-loop controls. 
\end{abstract}

\vspace{3mm}

 %60K35, 49L20
%60G40, 91A05,  49L20,  49L25.

\vspace{3mm}


\newpage

%\vspace{5mm}

 \section{Introduction}

\setcounter{equation}{0} \setcounter{Assumption}{0}
\setcounter{Theorem}{0} \setcounter{Proposition}{0}
\setcounter{Corollary}{0} \setcounter{Lemma}{0}
\setcounter{Definition}{0} \setcounter{Remark}{0}

The problem studied in this paper concerns the optimal control of mean-field stochastic differential equations (SDEs), also known as McKean-Vlasov equations.  This topic is closely related to the mean-field game  
(MFG) problem as originally formulated by Lasry and Lions in \cite{laslio07} and simultaneously by  Huang, Caines and Malham\'e in \cite{huaetal06}.  It aims at describing equilibrium states of large population of symmetric players  (particles) 
with mutual interactions of mean-field type, and we refer to \cite{cardel13} for a discussion pointing out the subtle differences between the notions of Nash equilibrium in MFG and Pareto optimality in the optimal control of McKean-Vlasov dynamics. 

While the analysis of McKean-Vlasov SDEs has a long history with the pioneering works by Kac \cite{kac56} and H. McKean \cite{mckean67}, and later on with papers in the  general framework of propagation of chaos, see e.g. 
\cite{sni89}, \cite{jouetal08},  the optimal control of McKean-Vlasov dynamics is a rather new problem, which attracts an increasing interest since the emergence of the MFG theory and its numerous applications in several areas outside physics, like economics and finance, biology, social interactions, networks. 
Actually,  it has been first studied in \cite{ahmdin01} by functional analysis method with a value function expressed in terms of the Nisio semigroup of operators.  More recently, se\-veral papers 
have adopted the  stochastic maximum (also called Pontryagin) principle for characterizing solutions to the controlled McKean-Vlasov systems in terms of an 
adjoint backward stochastic differential equation (BSDE) coupled with a forward SDE:  
see \cite{anddje10}, \cite{bucetal11}, \cite{yon13} with a state dynamics depending upon moments of the distribution, and \cite{cardel14} for a deep  investigation in a more general setting.  
Alternatively, and although the dynamics of mean-field SDEs  is non-Markovian, it is tempting  to use dynamic programming (DP) method (also called Bellman principle), 
which is known to be a powerful tool for standard Markovian stochastic control problem, see e.g. 
\cite{FleSon06}, \cite{Pha09},  and  does not require any convexity assumption usually imposed in Pontryagin principle. 
Indeed, mean-field type control problem was tackled by DP in \cite{laupir14} and \cite{benetal15} for specific  McKean-Vlasov SDE and cost functional, typically depending only upon statistics  
like its mean value or with uncontrolled diffusion coefficient, and  especially by assuming the existence at all times of a density for the marginal distribution of the state process. 
The key idea in both papers \cite{laupir14} and \cite{benetal15} is to reformulate the stochastic control problem 
with feedback strategy as a deterministic control problem involving the density of the marginal distribution, and then to derive a dynamic programming equation in the space of density functions.   
 
Inspired by the works \cite{benetal15} and \cite{laupir14},  the objective of this paper is to analyze in detail the dynamic programming method for the optimal control of mean-field SDEs where the drift, diffusion coeffi\-cients and running costs may depend both upon the joint  distribution of the state and of the control. This additional dependence related to the mean-field interaction on control is natural  in the context of McKean-Vlasov control problem, but has been few considered in the lite\-rature, see however \cite{yon13}  for a dependence only through the moments of the control. 
%Our paper can be viewed as the continuous time version of the discrete time problem stu\-died recently in \cite{phawei15}.    
By using closed-loop (also called feedback) controls, we first convert the stochastic optimal control problem  into a deterministic control problem where   the marginal distribution is the sole controlled  state variable, and we prove that dynamic programming holds in its general form.  The next step for exploiting the DP is to differentiate functions defined on the space  of probability measures. There are various notions of derivatives with respect to measures which have been developed in connection with the theory of optimal transport and using Wasserstein metric on the space of probability measures, see e.g. the monographs \cite{ambetal05}, \cite{vil09}.  For our purpose, we shall  use the notion of differentiability introduced by P.L. Lions in his lectures at the Coll\`ege de France \cite{lio12}, see also the helpful redacted notes \cite{car12}.  This notion of derivative is based on the lifting of functions defined on the space of square integrable proba\-bility measures into functions defined on the Hilbert space of square integrable random variables distributed according to the ``lifted" probability measure. It has been  used in \cite{cardel14} for differentiating the Hamiltonian function appearing in stochastic Pontryagin principle for controlled McKean-Vlasov dynamics. 
As usual in continuous time control problem, we  need a dynamic differential calculus for deriving the infinitesimal version of the DP, and shall rely on a special It\^o's chain rule for flows of probability measures as recently deve\-loped in \cite{buetal14} and \cite{chacridel15}, and used in \cite{cardel14b} 
for deriving the so-called Master equation in MFG.  We are then able to derive the dynamic progra\-mming Bellman equation for mean-field 
stochastic control problem. This infinite dimensional fully nonlinear partial differential equation (PDE) of second order in the Wassertein space of probability measures extends previous results in the literature \cite{benetal15}, \cite{cardel14b}, \cite{laupir14}:  it reduces in particular to the Bellman equation in the space of density functions  derived by Bensoussan, Frehse and Yam  \cite{benetal15b} when the marginal distribution admits a density, and  on the other hand, we notice that it differs from the Master equation for McKean-Vlasov control problem obtained by Carmona and Delarue in  \cite{cardel14b} where the value function is a function of both the state and its marginal distribution, and so with asso\-ciated PDE  in the state space comprising probability measures but also Euclidian vectors.  Following the traditional approach for stochastic control problem, we prove a verification theorem for the Bellman equation of the McKean-Vlasov control problem, which reduces to the classical 
Bellman equation in the case of  no mean-field interaction. We apply our verification theorem to the important class of linear quadratic (LQ) McKean-Vlasov control problems, addressed e.g. in \cite{yon13} and \cite{benetal16} by maximum principle and adjoint equations, and that we solve by a 
different approach where it turns out that derivations in the space of probability measures are quite tractable and lead to explicit classical solutions for  the Bellman equation.   We illustrate these results with two examples arising from finance: the mean-variance portfolio selection and an inter-bank systemic risk model, and retrieve the results obtained in \cite{lizho00}, \cite{fisliv15} and \cite{caretal14} by different methods. 

In general, there are  no classical solutions to the Bellman equation, and 
%and inspired by the various concepts of viscosity solutions for solutions to
%Hamilton-Jacobi in the Wasserstein space (and more generally in metric spaces), see \cite{ambetal05}, \cite{ganetal08}, \cite{car12}, \cite{ganswi15},  
we thus introduce a notion of viscosity solutions for the Bellman equation in the Wasserstein space of probability measures.  
There are several definitions of viscosity solutions for Hamilton Jacobi equations of first order in Wasserstein space and more gene\-rally in metric spaces, see e.g. \cite{ambetal05}, \cite{ganetal08}, 
\cite{fenkat09} or \cite{ganswi15}.  We adopt the approach in \cite{lio12}, and detailed in \cite{car12}, which consists, after the lifting identification between measures and random variables,  
in working in the Hilbert space of square integrable random variables  instead of working in  the Wasserstein space of probability measures, 
in order to use the various   tools developed for viscosity solutions in separable Hilbert spaces, in particular in our context,   for second order Hamilton-Jacobi equations, see \cite{lio88}, \cite{lio89b}, and the recent monograph \cite{fabgozswi15}.  We then prove the viscosity property of the value function and a comparison principle, hence uniqueness result,  for our  Bellman equation associated to  the McKean-Vlasov control problem.  
 
Finally, we consider the more general class of  open-loop controls instead of (Lipschitz) closed-loop  controls.  We derive the corresponding dynamic programming equation, and compare with the Bellman equation arising from McKean-Vlasov control problem with  feedback controls. 
  
The rest of the paper is organized as follows. Section 2  describes the McKean-Vlasov control problem and fix  the standing assumptions. 
In Section 3, we state  the dynamic programming principle  after the reformulation into a deterministic control problem, and derive the Bellman equation together with the proof of the verification theorem. We present in Section 4  the applications  to the LQ framework where explicit solutions are provided with two examples  arising from financial models.  Section 5  deals with viscosity solutions for the Bellman equation, and the last section considers the case of open-loop controls.

\section{McKean-Vlasov control problem}

\setcounter{equation}{0} \setcounter{Assumption}{0}
\setcounter{Theorem}{0} \setcounter{Proposition}{0}
\setcounter{Corollary}{0} \setcounter{Lemma}{0}
\setcounter{Definition}{0} \setcounter{Remark}{0}

Let us fix some probability space $(\Omega,\Fc,\P)$ on which is defined a $n$-dimensional Brownian motion $B$ $=$ $(B_t)_{0\leq t\leq T}$, and denote by $\F$ $=$ $(\Fc_t)_{0\leq t\leq T}$ its natural filtration, augmented with  an independent $\sigma$-algebra $\Fc_0$ $\subset$ $\Fc$. 
For each random variable $X$, we denote by $\P_{_X}$ its probability law (also called distribution)  under $\P$ (which is deterministic),
and by $\delta_X$  the  Dirac measure on $X$.
Given a normed space $(E,|.|)$, we denote by $\Pc_{_2}(E)$ the set of probability measures $\mu$ on $E$, which are square integrable, i.e.
$\|\mu\|_{_2}^2$ $:=$ $\int_E |x|^2 \mu(dx)$ $<$ $\infty$, and by $L^2(\Fc_0;E)$ ($=$ $L^2(\Omega,\Fc_0,\P;E)$) the set of square integrable random variables on $(\Omega,\Fc_0,\P)$.  
In the sequel, $E$ will be either $\R^d$, the state space, or $A$, the control space,  a subset of $\R^m$, or the product space $\R^d\times A$. 
We shall assume without loss of generality (see Remark \ref{remrich} below) 
that $\Fc_0$ is rich enough to carry  $E$-valued random variables with any arbitrary square integrable distribution, i.e. $\Pc_{_2}(E)$ $=$ $\{\P_\xi, \xi \in L^2(\Fc_0;E)\}$.  
%\marginpar{GIVE EXAMPLES FOR $(\Omega,\Fc_0,\P)$ to be rich enough} 

\begin{Remark} \label{remrich}
{\rm  A possible construction of a probability space, which is rich enough to satisfy the above conditions is the following. We consider a Polish space $\Omega_0$,  its Borel $\sigma$-algebra $\Fc_0$ and let $\P_0$ be  an atomless 
probability measure on $(\Omega_0,\Fc_0)$.  We consider another probability space $(\Omega_1,\Fc_1,\P_1)$ supporting a $n$-dimensional Brownian motion $B$ and denote by $\F^B$ $=$ $(\Fc_t^B)$ 
its natural filtration. By defining  $\Omega$ $=$ $\Omega_0\times\Omega_1$, $\Fc$ $=$ $\Fc_0\vee\Fc_1$,  $\P$ $=$ $\P_0\otimes\P_1$, and $\F$ $=$ $(\Fc_t)$ with $\Fc_t$ $=$ $\Fc_t^B \vee \Fc_0$, $0\leq t\leq T$,  
we then obtain that the filtered probability space $(\Omega,\Fc,\F,\P)$ satisfies the required  condition in the above framework. 
}
\ep
\end{Remark}

\vspace{1mm}

We also denote by $W_2$  the $2$-Wasserstein distance defined on $\Pc_{_2}(E)$ by
\beqs
W_2(\mu,\mu') &:=& \inf\Big\{ \Big( \int_{E\times E} |x-y|^2 \pi(dx,dy)\Big)^{1\over 2}:
\pi \in \Pc_{_2}(E\times E) \mbox{ with marginals } \mu \mbox{ and } \mu' \Big\} \\
&=& \inf\Big\{  \Big(\E|\xi-\xi'|^2\Big)^{1\over 2}: \;\; \xi,\xi' \in L^2(\Fc_0;E) \mbox{ with } \P_\xi = \mu, \; \P_{\xi'} = \mu' \Big\}. 
\enqs
%and we recall that $W_2(\P_\xi,\P_{\xi'})$ $\leq$ $\big(\E|\xi-\xi'|^2\big)^{1\over 2}$, for any $\xi$, $\xi'$ $\in$ $L^2(\Fc;E)$ ($=$ $L^2(\Omega,\Fc,\P;E)$). 

We consider a controlled stochastic dynamics of McKean-Vlasov type for the process $X^\alpha$ $=$ $(X_t^\alpha)_{0\leq t \leq T}$ valued in $\R^d$:
\beq \label{McKean}
dX_t^\alpha &=& b(t,X_t^\alpha,\alpha_t,\P_{_{(X_t^\alpha,\alpha_t)}}) dt +
\sigma(t,X_t^\alpha,\alpha_t,\P_{_{(X_t^\alpha,\alpha_t)}}) dB_t, \; \; X_0^\alpha \; = \; X_0,
\enq
where $X_0$ $\in$ $L^2(\Fc_0,\R^d)$, and  the control process $\alpha$ $=$ $(\alpha_t)_{0\leq t\leq T}$ is progressively measurable with values in a subset $A$ of $\R^m$, assumed for simplicity to contain the zero element.   
The coefficients $b$ and $\sigma$ are deterministic measurable  functions from $[0,T]\times\R^d\times A\times\Pc_{_2}(\R^d\times A)$ into $\R^d$ and $\R^{d\times n}$ respectively.  Notice here that the drift and diffusion coefficients
$b$, $\sigma$ of the controlled state process  do not depend only on the marginal distribution of the  state process $X_t$ at time $t$ but more generally on the joint  distribution of the state/control $(X_t,\alpha_t)$ at time $t$, which represents an additional
mean-field feature with respect to classical McKean-Vlasov equations.  We make the following assumption:

\vspace{2mm}

\hspace{-7mm} {\bf (H1)}
%we have
%\beqs
%\int_0^T |b(t,0,\delta_0,0,\delta_0)|^2  +  |\sigma(t,0,\delta_0,0,\delta_0)|^2 dt & < & \infty.
%\enqs
There exists some constant $C_{b,\sigma}$ $>$ $0$  s.t. for all $t$ $\in$ $[0,T]$, $x,x'$ $\in$ $\R^d$, $a,a'$ $\in$ $A$, 
$\lambda,\lambda'$ $\in$ $\Pc_{_2}(\R^d\times A)$,
\beqs
& & |b(t,x,a,\lambda) - b(t,x',a',\lambda')| + |\sigma(t,x,a,\lambda) - \sigma(t,x',a',\lambda')| \\
& \leq & C_{b,\sigma} \big[  |x-x'| + |a-a'| +  W_2(\lambda,\lambda') \big],
\enqs
and
\beqs
\int_0^T |b(t,0,0,\delta_{(0,0)})|^2  +  |\sigma(t,0,0,\delta_{(0,0)})|^2 dt  & < & \infty.
\enqs

\vspace{2mm}

Condition {\bf(H1)} ensures that for any control process $\alpha$, which is square integrable, i.e. $\E[\int_0^T |\alpha_t|^2 dt]$ $<$ $\infty$,  there exists a unique solution
$X^\alpha$ to \reff{McKean}, and moreover this solution satisfies (see e.g. \cite{sni89} or \cite{jouetal08}):
\beq \label{momentX}
\E \big[ \sup_{0\leq t\leq T} |X_t^\alpha|^2  \big] & \leq  & C\Big( 1 + \E|X_0|^2 + \E\big[\int_0^T |\alpha_t|^2 dt\big] \Big)  \; < \; \infty. 
\enq
In the sequel of the paper, we stress the dependence of $X^\alpha$ on $\alpha$ if needed, but most often,
we shall omit  this  dependence and simply write  $X$ $=$ $X^\alpha$   when there is no ambiguity.

\vspace{2mm}

The  cost functional associated to the McKean-Vlasov equation \reff{McKean} is
\beq \label{gainJ}
J(\alpha) &:=& \E \Big[ \int_0^T f(t,X_t,\alpha_t,\P_{_{(X_t,\alpha_t)}}) dt + g(X_T,\P_{_{X_T}}) \Big]
\enq
for a square integrable control process $\alpha$.  The running cost function $f$ is a deterministic real-valued function on
$[0,T]\times\R^d\times A\times\Pc_{_2}(\R^d\times A)$ and the terminal gain function $g$ is a deterministic real-valued function on
$\R^d\times\Pc_{_2}(\R^d)$.  We shall assume the following quadratic condition on $f$, $g$:

\vspace{2mm}

\hspace{-7mm} {\bf (H2)} There exists  some constant $C_{f,g}$ $>$ $0$ s.t. for all $t$ $\in$ $[0,T]$, $x$ $\in$ $\R^d$, $a$ $\in$ $A$,
$\mu$ $\in$ $\Pc_{_2}(\R^d)$, $\lambda$ $\in$ $\Pc_{_2}(\R^d\times A)$,
\beqs
|f(t,x,a,\lambda)| + |g(x,\mu)| & \leq & C_{f,g} \big(1 + |x|^2 + |a|^2 + \|\mu\|^2_{_2} + \|\lambda\|^2_{_2} \big).
\enqs

\vspace{2mm}

Under Condition {\bf (H2)}, and from  \reff{momentX}, we see that $J(\alpha)$ is well-defined and finite for any square integrable control process $\alpha$.
The stochastic control problem of interest in this paper is to minimize  the cost functional:
\beq \label{defV0}
V_0 &:=& \inf_{\alpha\in\Ac} J(\alpha),
\enq
over a set of admissible controls $\Ac$ to be precised later. 
%Problem \reff{defV0} arises in the asymptotic study of collective behaviors  of a large number of players (particles) resulting from mean-field interactions, %and when a center decides of the same control policy  for all symmetric players. We refer to \cite{cardel13} for a detailed discussion.
%In the case of drift and diffusion coefficients not depending on the marginal distribution of the control
%$\alpha$, problem \reff{defV0} was studied in \cite{cardel14}  by stochastic maximum principle method leading to  a system of (coupled)  forward %backward stochastic differential equations for the solution to the control problem

%On the other hand, dynamic programming is a well-known methodology  for studying classical Markovian stochastic control problem, see e.g.
%\cite{FleSon06} or \cite{Pha09}, and  our purpose in this paper is to  develop this approach in the general McKean-Vlasov framework, and show how it %leads to a Bellman equation, which can be explicitly solved in some particular cases.

\vspace{2mm}

\noindent {\bf Notations}:  We denote by $x.y$ the scalar product of two Euclidian vectors $x$ and $y$, and by $M\trans$ the transpose of a matrix or vector $M$.  For any $\mu$ $\in$ $\Pc_{_2}(E)$, $F$ Euclidian space, we denote 
by $L^2_\mu(F)$ the set of measurable functions $\varphi$ $:$ $E$ $\rightarrow$ $F$ which are square integrable with respect to $\mu$,  and we set
\beqs
<\varphi,\mu> &:=& \int_{E} \varphi(x) \mu(dx). 
\enqs
We also denote by $L^\infty_\mu(F)$ the set of measurable functions $\varphi$ $:$ $E$ $\rightarrow$ $F$ which are bounded $\mu$ a.e., and  
$\|\varphi\|_{_ {\infty}}$ denotes the essential supremum of $\varphi$ $\in$ $L_\mu^\infty(F)$.

\section{Dynamic programming and Bellman equation}

\setcounter{equation}{0} \setcounter{Assumption}{0}
\setcounter{Theorem}{0} \setcounter{Proposition}{0}
\setcounter{Corollary}{0} \setcounter{Lemma}{0}
\setcounter{Definition}{0} \setcounter{Remark}{0}

\subsection{Dynamic programming principle}

In this paragraph, we make the standing assumptions {\bf (H1)}-{\bf (H2)}, and our purpose is to show  that dynamic programming principle holds
for problem \reff{defV0},  which we would like to combine with some Markov property of the controlled state process. However, notice that  the McKean-Vlasov type dependence on the dynamics of the state process rules out the standard Markov property  of the controlled process $(X_t)_t$. Actually, this Markov property can be restored by considering its  probability law  $(\P_{_{X_t}})_t$.  To be more precise and for the sake of definiteness, we shall restrict ourselves to controls
$\alpha$ $=$ $(\alpha_t)_{0\leq t\leq T}$ given in closed loop (or feedback) form:
\beq \label{closed}
\alpha_t &=& \tilde\alpha(t,X_t,\P_{_{X_t}}), \;\;\;\;\;  0 \leq t \leq T,
\enq
for some deterministic measurable function $\tilde\alpha(t,x,\mu)$ defined on $[0,T]\times\R^d\times\Pc_{_2}(\R^d)$. 
We shall discuss in the last section how one deal more generally with  open-loop controls. 
We  denote by $Lip([0,T]\times\R^d\times\Pc_{_2}(\R^d);A)$
the set of deterministic measurable functions $\tilde\alpha$ on  $[0,T]\times\R^d\times\Pc_{_2}(\R^d)$, valued in $A$, which are Lipschitz in $(x,\mu)$,
and satisfy a linear growth condition on $(x,\mu)$, uniformly on $t$ $\in$ $[0,T]$, i.e. there exists some positive constant $C_{\tilde\alpha}$ s.t. for all $t$ $\in$ $[0,T]$, $x,x'$ $\in$
$\R^d$, $\mu,\mu'$ $\in$ $\Pc_{_2}(\R^d)$,
\beqs
|\tilde\alpha(t,x,\mu) - \tilde\alpha(t,x',\mu') | & \leq & C_{\tilde\alpha} \big( | x- x'| + W_2(\mu,\mu') \big), \\
\int_0^T |\tilde\alpha(t,0,\delta_0)|^2 dt  & < & \infty.
\enqs
Notice that for any $\tilde\alpha$ $\in$ $Lip([0,T]\times\R^d\times\Pc_{_2}(\R^d);A)$,  and under the Lipschitz condition in {\bf (H1)},  there exists a unique solution to the SDE:
\beq
dX_t &=& b(t,X_t,\tilde\alpha(t,X_t,\P_{_{X_t}}),\P_{_{(X_t,\tilde\alpha(t,X_t,\P_{_{_{X_t}}}))}}) dt  \nonumber \\
& & \; \;\; + \; \sigma(t,X_t,\tilde\alpha(t,X_t,\P_{_{X_t}}),\P_{_{(X_t,\tilde\alpha(t,X_t,\P_{_{X_t}}))}}) dB_t,   \label{Xfeedback}
\enq
starting from some square integrable random variable,  and this solution satisfies the square integrability condition \reff{momentX}. The set $\Ac$ of so-called admissible controls $\alpha$ is then defined as the set of control processes
$\alpha$ of feedback form \reff{closed} with $\tilde\alpha$ $\in$ $Lip([0,T]\times\R^d\times\Pc_{_2}(\R^d);A)$. We shall often identify $\alpha$ $\in$ $\Ac$ with $\tilde\alpha$ in
$Lip([0,T]\times\R^d\times\Pc_{_2}(\R^d);A)$ via \reff{closed}, and we see that any $\alpha$ in $\Ac$ is square-integrable: $\E[\int_0^T |\alpha_t|^2 dt]$ $<$ $\infty$, by \reff{momentX} and
Gronwall's lemma.

Let us now check  the flow property of the marginal distribution process $\P_{_{X_t}}$ $=$ $\P_{_{X_t^\alpha}}$ for any admissible control $\alpha$ in 
$\Ac$. For any $\tilde\alpha$ $\in$ $L(\R^d;A)$, the set of Lipschitz functions from $\R^d$ into $A$, we denote by 
$Id\tilde\alpha$ the function 
\beqs
Id\tilde\alpha : \R^d & \rightarrow & \R^d\times A \\
x & \mapsto &   (x,\tilde\alpha(x)).
\enqs
We observe that the joint distribution $\P_{_{(X_t,\alpha_t)}}$ associated to a feedback control $\alpha$ $\in$ $\Ac$  is equal to 
the image by $Id\tilde\alpha(t,.,\P_{_{X_t}})$ of the marginal distribution $\P_{_{X_t}}$ of the controlled  state process $X$, i.e. 
$\P_{_{(X_t,\alpha_t)}}$ $=$ $Id\tilde\alpha(t,.,\P_{_{X_t}})\star\P_{_{X_t}}$, where $\star$ denotes the standard pushforward of measures: 
for any $\tilde\alpha$ $\in$ $L(\R^d;A)$, and $\mu$  $\in$ $\Pc_{_2}(\R^d)$:
\beqs
(Id\tilde\alpha\star\mu)(B) &=& \mu\big(Id\tilde\alpha^{-1}(B) \big), \;\;\; \forall B \in \Bc(\R^d\times A).
\enqs
We  consider the dynamic version of \reff{Xfeedback} starting at time $t$ $\in$ $[0,T]$ from $\xi$ $\in$ $L^2(\Fc_t;\R^d)$, which is then written as:
%(set of square-integrable $\Fc_t$-measurable random variables valued in $\R^d$):
\beq
X_s^{t,\xi}  &=& \xi  \; + \;   \int_t^s b(r,X_r^{t,\xi},\tilde\alpha(r,X_r^{t,\xi},\P_{_{X_r^{t,\xi}}}), Id\tilde\alpha(r,.,\P_{_{X_r^{t,\xi}}})\star\P_{_{X_r^{t,\xi}}}) dr
\label{Xfeedbackxi} \\
& & \; + \;  \int_t^s \sigma(r,X_r^{t,\xi},\tilde\alpha(r,X_r^{t,\xi},\P_{_{X_r^{t,\xi}}}),Id\tilde\alpha(r,.,\P_{_{X_r^{t,\xi}}})\star\P_{_{X_r^{t,\xi}}}) dB_r, \;\;\; t \leq s \leq T.   \nonumber
\enq
Existence and uniqueness of  a solution to \reff{Xfeedbackxi} implies the flow property:
\beq \label{flowX}
X_s^{t,\xi} &=& X_s^{\theta,X_\theta^{t,\xi}}, \;\;\;\;\;    \forall \;  0 \leq  t \leq \theta \leq s   \leq T,  \; \xi \in L^2(\Fc_t;\R^d).
\enq
Moreover,  as pointed out in Remark 3.1 in \cite{buetal14} (see also the remark following (2.3) in \cite{chacridel15}),  the solution to \reff{Xfeedbackxi} is also unique in law from which it follows that
the law of $X^{t,\xi}$  depends on $\xi$ only through its law $\P_{_\xi}$.  Therefore, we can define
\beq \label{defPmu}
\P_s^{t,\mu} &:=& \P_{_{X_s^{t,\xi}}}, \;\;\; \mbox{ for }  \;\;\;   0 \leq t\leq  s\leq T, \;\; \mu \; = \;  \P_{_{\xi}} \; \in \; \Pc_{_2}(\R^d),
\enq
As a consequence of the flow property \reff{flowX}, and recalling that $\Pc_{_2}(\R^d)$ $=$ $\{\P_\xi, \xi \in L^2(\Fc_0;\R^d)\}$,
it is clear that we also get the flow property for the marginal distribution process:
\beq \label{flowP}
\P_s^{t,\mu} &=& \P_s^{\theta,\P_\theta^{t,\mu}},  \;\;\;\;\;    \forall \;  0 \leq  t \leq \theta \leq s   \leq T, \; \mu \in \Pc_{_2}(\R^d).
\enq
Recall that the process $X^{t,\xi}$, hence also the law process $\P^{t,\mu}$ depends on the feedback control $\alpha$ $\in$ $\Ac$, and if needed, we shall stress the dependence on $\alpha$ by writing $\P^{t,\mu,\alpha}$.

We  next  show that the initial stochastic control problem can be reduced to a deter\-ministic control problem. Indeed,  by definition of the marginal distribution $\P_{_{X_t}}$, recalling that $\P_{_{(X_t,\alpha_t)}}$ $=$ $Id\tilde\alpha(t,.,\P_{_{X_t}})\star\P_{_{X_t}}$, and Fubini's theorem, we see  that the cost functional can be written for any admissible control $\alpha$ $\in$ $\Ac$ as:
\beqs
J(\alpha) &=& \int_0^T \hat f(t,\P_{_{X_t}},\tilde\alpha(t,.,\P_{_{X_t}})) dt + \hat g(\P_{_{X_T}}),
\enqs
where the function $\hat f$ is defined on $[0,T]\times\Pc_{_2}(\R^d)\times L(\R^d;A)$ and $\hat g$ is defined on $\Pc_{_2}(\R^d)$ by
\beq \label{defhatfg}
\hat f(t,\mu,\tilde\alpha) \; :=\;  <f(t,.,\tilde\alpha(.),Id\tilde\alpha\star\mu),\mu >,  & &
%\int_{\R^d} f(x,\mu,\tilde\alpha(x),\tilde\alpha\star\mu) \mu(dx),
\;\;\; \hat g(\mu) \; := \; <g(.,\mu),\mu >.
%\int_{\R^d} g(x,\mu) \mu(dx).
\enq
%by using the  notation $<\varphi,\mu>$ $=$ $\int_{\R^d} \varphi(x) \mu(dx)$ for any
%$\mu$ $\in$ $\Pc_{_2}(\R^d)$, and $\varphi$ measurable function with quadratic growth condition.
We have thus transformed the initial control problem \reff{defV0} into a deterministic control problem involving the infinite dimensional controlled
marginal distribution process valued in $\Pc_{_2}(\R^d)$. In view of the flow property \reff{flowP}, it is then natural to define the value function
\beq \label{defv}
v(t,\mu) & := & \inf_{\alpha\in\Ac} \Big[ \int_t^T \hat f(s,\P_s^{t,\mu},\tilde\alpha(s,.,\P_s^{t,\mu})) ds + \hat g(\P_T^{t,\mu}) \Big], \;\; t \in [0,T], \; \mu \in \Pc_{_2}(\R^d),
\enq
so that the initial control problem in \reff{defV0} is given by: $V_0$ $=$ $v(0,\P_{_{X_0}})$.   It is clear that $v(t,\mu)$ $<$ $\infty$, and we shall assume that
\beq \label{Vfini}
v(t,\mu)   & > & - \infty, \;\;\;\;\;  \forall \; t \in [0,T], \; \mu \in \Pc_{_2}(\R^d).
\enq

\begin{Remark} \label{remvfini}
{\rm The finiteness condition \reff{Vfini} can be checked a priori directly from the assump\-tions on the model. For example, when $f$, $g$, hence $\hat f$, $\hat g$,
are lower-bounded functions, condition  \reff{Vfini} clearly holds.  Another  example is the case when $f(t,x,a,\lambda)$, and $g(x,\mu)$ are lower bounded by a quadratic function  in $x$, $\mu$, and $\lambda$ (uniformly in $(t,a)$) so that
\beqs
\hat f(t,\mu,\tilde\alpha) + \hat g(x,\mu) & \geq & - C \big( 1 + \|\mu\|_{_2} \big), \;\;\; \forall \mu \in \Pc_{_2}(\R^d), \;  \tilde\alpha\in L(\R^d;A),
\enqs
and we are able to derive moment estimates on the controlled process $X$, uniformly in $\alpha$: $\big\|\P_{_s}^{t,\mu}\big\|_{_2}^2$ $=$ $\E[|X_s^{t,\xi}|^2]$ $\leq$ $C(1+\|\mu\|_{_2}^2)$, (for $\mu$ $=$ $\P_{_\xi}$)
which arises typically from \reff{momentX} when $A$ is bounded. Then, it is clear that \reff{Vfini} holds true.
Otherwise, this finiteness  condition can be checked a posteriori from a verification theorem, see Theorem \ref{theoverif}.
}
\ep
\end{Remark}

The  dynamic programming principle (DPP) for the  deterministic control pro\-blem \reff{defv} takes the following formulation:

\begin{Theorem} \label{theoDPP}
(Dynamic Programming Principle)

\noindent Under \reff{Vfini}, we have for all  $0\leq t\leq \theta\leq T$, $\mu$ $\in$ $\Pc_{_2}(\R^d)$:
\beq \label{DPPdeter}
v(t,\mu)  &=& \Inf_{\alpha\in\Ac } \Big[ \int_t^\theta \hat f(s,\P_s^{t,\mu},\tilde\alpha(s,.,\P_s^{t,\mu}))  ds \; + \;  v(\theta,\P_\theta^{t,\mu}) \Big].
\enq
\end{Theorem}
{\bf Proof.} In the context of deterministic control problem, the proof of the DPP is elementary and does not require any measurable selection arguments. For sake of completeness, we provide  it.  Denote by $J(t,\mu,\alpha)$ the cost functional:
\beqs
J(t,\mu,\alpha) &:=& \int_t^T \hat f(s,\P_{s}^{t,\mu,\alpha},\tilde\alpha(s,.,\P_{s}^{t,\mu,\alpha}))ds + \hat g(\P_{T}^{t,\mu,\alpha}), \;\; 0 \leq t\leq T, \; \mu \in \Pc_{_2}(\R^d), \alpha\in \Ac,
\enqs
so that $v(t,\mu)$ $=$ $\inf_{\alpha\in\Ac}J(t,\mu,\alpha)$, and by $w(t,\mu)$ the r.h.s. of \reff{DPPdeter} (here we stress the dependence of
the controlled marginal distribution process $\P^{t,\mu,\alpha}$ on $\alpha$).   Then,
\beqs
w(t,\mu) & =&  \Inf_{\alpha\in\Ac} \big[ \int_t^\theta  \hat f(s,\P_{s}^{t,\mu,\alpha},\tilde\alpha(s,.,\P_{s}^{t,\mu,\alpha})) ds
+  \inf_{\beta\in\Ac} J(\theta,P_\theta^{t,\mu,\alpha},\beta) \big] \\
&=& \Inf_{\alpha\in\Ac} \inf_{\beta\in\Ac} \big[  \int_t^\theta  \hat f(s,\P_{s}^{t,\mu,\alpha},\tilde\alpha(s,.,\P_{s}^{t,\mu,\alpha})) ds +   J(\theta,P_\theta^{t,\mu,\alpha},\beta) \big] \\
&=&  \Inf_{\alpha\in\Ac} \inf_{\beta\in\Ac} \big[  \int_t^\theta \hat f(s,\P_{s}^{t,\mu,\gamma[\alpha,\beta]},\tilde\gamma[\alpha,\beta](s,.,\P_{s}^{t,\mu,\gamma[\alpha,\beta]})) ds
+   J(\theta,P_\theta^{t,\mu,\gamma[\alpha,\beta]},\gamma[\alpha,\beta]) \big]
%&=& \inf_{\gamma\in  \{ \Ac_{\theta}(\beta): \beta \in \Ac_t(\alpha)\}   }  J_{t}(\gamma),
\enqs
where we define  $\gamma[\alpha,\beta]$ $\in$ $\Ac$ by:  $\tilde\gamma[\alpha,\beta](s,.)$ $=$ $\tilde\alpha(s,.) 1_{0\leq s\leq\theta} + \tilde\beta(s,.) 1_{\theta<s\leq T}$. Now, it is clear that when $\alpha, \beta$ run over $\Ac$, then
$\gamma[\alpha,\beta]$ also runs over $\Ac$, and so:
\beqs
w(t,\mu) & =& \inf_{\gamma\in\Ac}  \big[  \int_t^\theta \hat f(s,\P_{s}^{t,\mu,\gamma},\tilde\gamma(s,.,\P_{s}^{t,\mu,\gamma})) ds +   J(\theta,P_\theta^{t,\mu,\gamma},\gamma) \big] \\
&=&  \inf_{\gamma\in\Ac}  \big[  \int_t^\theta \hat f(s,\P_{s}^{t,\mu},\tilde\gamma(s,.,\P_{s}^{t,\mu})) ds + \int_\theta^T \hat f(s,\P_{s}^{\theta,\P_\theta^{t,\mu}},\tilde\gamma(s,.,\P_{s}^{\theta,\P_\theta^{t,\mu}})) + \hat g(\P_{T}^{\theta,\P_\theta^{t,\mu}})
\big]  \\
&=&   \inf_{\gamma\in\Ac}  \big[  \int_t^\theta \hat f(s,\P_{s}^{t,\mu},\tilde\gamma(s,.,\P_{s}^{t,\mu})) ds + \int_\theta^T \hat f(s,\P_{s}^{t,\mu},\tilde\gamma(s,.,\P_{s}^{t,\mu})) + \hat g(\P_{T}^{t,\mu}) \big],
\enqs
by the flow property \reff{flowP} (here we have omitted in the second and third line the dependence of $\P_s$ in $\gamma$).   This proves the required equality: $w(t,\mu)$ $=$ $v(t,\mu)$.
\ep

\vspace{3mm}

\begin{Remark} \label{remoptcontrol}
{\rm Problem \reff{defV0}  includes the case where the cost functional in \reff{gainJ} is a nonlinear function of the expected value of the state process, i.e.  the running cost functions and the terminal gain function are in the form:
$f(t,X_t,\alpha_t,\P_{_{(X_t,\alpha_t)}})$ $=$ $\bar f(t,X_t,\E[X_t],\alpha_t)$, $t$ $\in$ $[0,T]$,  $g(X_T,\P_{_{X_T}})$ $=$ $\bar g(X_T,\E[X_T])$,  which arises for example in mean-variance problem (see Section \ref{secappli}).  It is claimed in \cite{bjomur08} and \cite{yon13}  that Bellman optimality principle does not hold, and therefore the problem is time-inconsistent.  This is correct when one takes into account only the
state process $X$ (that is its realization), since it is not Markovian, but as shown in this section, dynamic programming principle holds true whenever we consider the marginal distribution as state variable.  This gives more information and the price to paid is the infinite-dimensional feature of the marginal distribution state variable.
}
\ep
\end{Remark}

\subsection{Bellman equation}

The purpose of this paragraph  is to  derive from the dynamic programming principle \reff{DPPdeter},  a partial differential equation (PDE) for the value function $v(t,\mu)$, called Bellman equation. We shall rely on the notion of derivative with respect to a probability measure, as introduced by P.L. Lions in his course at Coll\`ege de France, and detailed   in the lecture notes  \cite{car12}.

This notion is based on the lifting of functions $u$ $:$ $\Pc_{_2}(\R^d)$ $\rightarrow$ $\R$ into functions $U$ defined on $L^2(\Fc_0;\R^d)$ by $U(X)$ $=$ $u(\P_{_X})$. We say that $u$ is differentiable (resp. $\Cc^1$) on
$\Pc_{_2}(\R^d)$ if the lift $U$ is Fr\'echet differentiable (resp. Fr\'echet differentiable with continuous derivatives) on  $L^2(\Fc_0;\R^d)$. In this case, the Fr\'echet derivative $[DU](X)$, viewed as an element $DU(X)$ of $
L^2(\Fc_0;\R^d)$  by Riesz' theorem:  $[DU](X)(Y)$ $=$ $\E[DU(X).Y]$,  can be represented as
\beq \label{Uu1}
DU(X) &=& \partial_\mu u(\P_{_X})(X),
\enq
for some function  $\partial_\mu u(\P_{_X})$ $:$ $\R^d$ $\rightarrow$ $\R^d$,  which is  called derivative of $u$ at $\mu$ $=$ $\P_{_X}$.  Moreover, $\partial_\mu u(\mu)$ $\in$ $L^2_\mu(\R^d)$ for 
$\mu$ $\in$ $\Pc_{_2}(\R^d)$ $=$ $\{\P_{_X}, X \in L^2(\Fc_0;\R^d)\}$. 
Following \cite{chacridel15}, we say that $u$ is partially $\Cc^2$ if it is $\Cc^1$, and one can find, for any $\mu$ $\in$ $\Pc_{_2}(\R^d)$, 
a continuous version of the mapping $x\in\R^d$ $\mapsto$ $\partial_\mu u(\mu)(x)$, such that the mapping
$(\mu,x)$ $\in$ $\Pc_{_2}(\R^d)\times\R^d$ $\mapsto$ $\partial_\mu u(\mu)(x)$  is continuous at any point $(\mu,x)$ such that $x$ $\in$ Supp$(\mu)$, and if for any 
$\mu$ $\in$ $\Pc_{_2}(\R^d)$, the mapping
$x$ $\in$ $\R^d$ $\mapsto$  $\partial_\mu u(\mu)(x)$ is differentiable, its derivative being jointly continuous at any point  $(\mu,x)$ such that $x$ $\in$ Supp$(\mu)$. The gradient is then denoted by $\partial_x  \partial_\mu u(\mu)(x)$  $\in$ $\S^{d}$, the set of symmetric matrices in $\R^{d\times d}$.  
% and  $\partial_x  \partial_\mu u(\mu)$ $\in$ $L_\mu^2(\R^{d\times d})$. 
We say that $u$ $\in$ $\Cc^2_b(\Pc_{_2}(\R^d))$ if it is partially $\Cc^2$,  $\partial_x  \partial_\mu u(\mu)$ $\in$ $L_\mu^\infty(\S^{d})$, 
and for any compact set $\Kc$ of $\Pc_{_2}(\R^d)$, we have
%there exists some constant $C_u$ s.t. for all $\mu$ $\in$ $\Pc_{_2}(\R^d)$,
\beqs
 \sup_{ \mu \in \Kc } \Big[ \int_{\R^d} \big| \partial_\mu u(\mu)(x) |^2\mu(dx)  +
%\int_{\R^d}  \big| \partial_x \partial_\mu u(\mu)(x) |^2 \mu(dx)  
\big \| \partial_x \partial_\mu u(\mu)\|_{_\infty}
\Big]  & < & \infty.
%& \leq & C_u\big( 1 + \|\mu\|^2_{_2} \big).
\enqs
%In this case, notice that $\partial_\mu u(\mu)$ (resp.  $\partial_x \partial_\mu u(\mu)$)  lies in $L_\mu^2(\R^d;\R^d)$ (resp. $L_\mu^2(\R^d;\R^{d\times d})$) the set functions from $\R^d$ into $\R^d$ (resp. $\R^{d\times d}$),  
%and square-integrable w.r.t. $\mu$.  
As shown in \cite{chacridel15}, if the lifted function $U$ is twice continuously Fr\'echet differentiable on   $L^2(\Fc_0;\R^d)$ with Lipschitz Fr\'echet derivative, 
then $u$ lies in $\Cc^2_b(\Pc_{_2}(\R^d))$.  In this case, the second Fr\'echet derivative $D^2 U(X)$  is identified indifferently by Riesz' theorem  as a  bilinear form on $L^2(\Fc_0;\R^d)$ or  as a symmetric operator  (hence bounded) on  $L^2(\Fc_0;\R^d)$, denoted  by $D^2U(X)$ $\in$ $S(L^2(\Fc_0;\R^d))$,  and we have the relation (see Appendix A.2  in \cite{cardel14b}): 
\beq \label{Uu2}
 \E \Big[ D^2U(X)(YN) . YN \Big]   &=& \E\Big[ {\rm tr} \big( \partial_x \partial_\mu u(\P_{_X})(X) YY\trans \big) \Big],
\enq
for any $X$ $\in$ $L^2(\Fc_0;\R^d)$, $Y$ $\in$ $L^2(\Fc_0;\R^{d\times q})$,  and where $N$ $\in$ $L^2(\Fc_0;\R^q)$ is  independent of $(X,Y)$ with zero mean 
and unit variance. 
%\marginpar{CHECK IF $C^2$ Frechet differentiable implies $C^2$ in the space of measures}  

We shall need a chain rule (or It\^o's formula)  for functions defined on $\Pc_{_2}(\R^d)$, proved independently in \cite{buetal14} and \cite{chacridel15}, see also the Appendix in \cite{cardel14b},  and that we recall here.  Let us consider an $\R^d$-valued It\^o process
\beqs
dX_t &=& b_t dt + \sigma_t dB_t, \;\;\; X_0 \in L^2(\Fc_0;\R^d),
\enqs
where $(b_t)$ and $(\sigma_t)$ are progressively measurable processes with respect to the filtration generated by the $n$-dimensional Brownian motion $B$,  valued respectively  in $\R^d$ and $\R^{d\times n}$, and satisfying the integrability condition:
\beq \label{condsigma}
\E \Big[ \int_0^T |b_t|^2 + |\sigma_t|^2 dt \Big] & < & \infty.
\enq
Let $u$ $\in$  $\Cc^2_b(\Pc_{_2}(\R^d))$. Then, for all $t$ $\in$ $[0,T]$,
\beq \label{Ito}
u(\P_{_{X_t}}) &=& u(\P_{_{X_0}})  \; + \; \int_0^t  \E \big[ \partial_\mu u(\P_{_{X_s}})(X_s).b_s     +  \frac{1}{2}   {\rm tr}\big(\partial_x\partial_\mu u(\P_{_{X_s}})(X_s) \sigma_s\sigma_s\trans\big) \big]  ds. 
\enq
%where we denote by $.$ the scalar product in $\R^d$ and by $\sigma\trans$ the transpose of a matrix $\sigma$.

We have now the ingredients for deriving the Bellman equation associated to the DPP \reff{DPPdeter}, and it turns out that it  takes   the following form:
\begin{equation} \label{HJB}
\left\{
\begin{array}{rcl}
\partial_t v \; + \;   \Inf_{\tilde\alpha\in L(\R^d;A)} \Big[ \hat f(t,\mu,\tilde\alpha)  \; + \;   <\Lc_t^{\tilde\alpha} v(t,\mu) , \mu> \Big]  &=& 0, \; \mbox{ on } \; [0,T)\times\Pc_{_2}(\R^d),    \\
v(T,.) &=& \hat g,  \; \mbox{ on } \; \Pc_{_2}(\R^d)
\end{array}
\right.
\end{equation}
where  for $\tilde\alpha$ $\in$ $L(\R^d;A)$, $\varphi$ $\in$ $\Cc^2_b(\Pc_{_2}(\R^d))$ and   $(t,\mu)$ $\in$ $[0,T]\times\Pc_{_2}(\R^d)$,
$\Lc_t^{\tilde\alpha}\varphi(\mu)$ $\in$ $L^2_\mu(\R)$  is the function: $\R^d$ $\rightarrow$ $\R$  defined by
\beq
\Lc_t^{\tilde\alpha} \varphi(\mu)  (x)  &:=&  \partial_\mu \varphi(\mu)(x).b(t,x,\tilde\alpha(x),Id\tilde\alpha\star\mu)   \nonumber \\
& & \;\;\; + \;   \frac{1}{2}   {\rm tr}\big(\partial_x\partial_\mu \varphi(\mu)(x) \sigma\sigma\trans (t,x,\tilde\alpha(x),Id\tilde\alpha\star\mu) \big). \label{Lcalpha}
\enq

%\subsection{Verification theorem}

In the spirit of classical verification theorem for stochastic control of diffusion processes, we  prove  the following result in our McKean-Vlasov control framework, which 
is a consequence of  It\^o's formula for functions defined on the Wasserstein space. 
%shall rely as usual  on It\^o's formula.  
%Actually, in order to satisfy the integrability conditions \reff{condsigma}, we need to strengthen the growth condition on $\sigma$ by requiring the existence of some positive constant $C_\sigma$ s.t.  for all
%$(t,x,\mu,a,\lambda)$ $\in$ $[0,T]\times\R^d\times\Pc_{_2}(\R^d)\times A\times\Pc_{_2}(A)$:
%\beq \label{sigmabor}
%|\sigma(t,x,\mu,a,\lambda)|^2  & \leq & C_\sigma\big(1 + |x| + |a| +  \|\mu\|_{_2} +  \|\lambda\|_{_2} \big).
%\enq

\begin{Proposition} \label{theoverif} (Verification theorem)

\noindent Let $w$ $:$ $[0,T]\times\Pc_{_2}(\R^d)$ $\rightarrow$ $\R$ be a function in $\Cc_b^{1,2}([0,T]\times\Pc_{_2}(\R^d))$, i.e.   $w$ is continuous on $[0,T]\times\Pc_{_2}(\R^d)$,  $w(t,.)$ $\in$ $\Cc_b^2(\Pc_{_2}(\R^d))$, for all $t$ $\in$ $[0,T]$, and $w(.,\mu)$ $\in$ $\Cc^1([0,T))$.
Suppose that  $w$ is solution to \reff{HJB}, and there exists  for all $(t,\mu)$ $\in$ $[0,T)\times\Pc_{_2}(\R^d)$ an element $\tilde\alpha^*(t,.,\mu)$ $\in$ $L(\R^d;A)$ attaining the infimum in \reff{HJB} s.t.
the mapping $(t,x,\mu)$ $\mapsto$ $\tilde\alpha^*(t,x,\mu)$ $\in$  $Lip([0,T]\times\R^d\times\Pc_{_2}(\R^d);A)$.
Then,  $w$ $=$ $v$, and the feedback control $\alpha^*$ $\in$ $\Ac$ defined by
\beqs
\alpha_t^* &=& \tilde\alpha^*(t,X_t,\P_{_{X_t}}), \;\;\; 0 \leq t <  T,
\enqs
is an optimal control, i.e.  $V_0$ $=$ $J(\alpha^*)$.
\end{Proposition}
{\bf Proof.}   Fix $(t,\mu=\P_\xi)$ $\in$ $[0,T)\times\Pc_{_2}(\R^d)$, and consider some arbitrary feedback control $\alpha$ $\in$ $\Ac$
associated to  $X^{t,\xi}$ the solution to the controlled SDE \reff{Xfeedbackxi}.  Under condition {\bf (H1)}, we have the standard estimate
\beqs
\E \big[ \sup_{t\leq s\leq T} |X_s^{t,\xi}|^2  \big] & \leq  & C\big( 1 + \E|\xi|^2 \big)  \; < \; \infty,
\enqs
which implies  that
%under \reff{sigmabor} that
\beqs
& & \E \Big[ \int_t^T
\big|b(s,X_s^{t,\xi},\tilde\alpha(s,X_s^{t,\xi},\P_{_{X_s^{t,\xi}}}),Id\tilde\alpha(s,.,\P_{_{X_s^{t,\xi}}})\star\P_{_{X_s^{t,\xi}}})\big|^2 \\
& & \;\;\;\;\;\;\; + \; \big|\sigma(s,X_s^{t,\xi},\tilde\alpha(s,X_s^{t,\xi},\P_{_{X_s^{t,\xi}}}),Id\tilde\alpha(s,.,\P_{_{X_s^{t,\xi}}})\star\P_{_{X_s^{t,\xi}}})\big|^2 ds \Big]  \; \; <  \; \; \infty.
\enqs
One can then apply  the It\^o's formula  \reff{Ito} to
$w(s,\P_{_{X_s^{t,\xi}}})$ $=$ $w(s,\P_s^{t,\mu})$  (with the definition \reff{defPmu}) between $s$ $=$ $t$ and $s$ $=$ $T$, and obtain
\beq
w(T,\P_T^{t,\mu}) &=& w(t,\mu) \; + \; \int_t^T \Dt{w}(s,\P_s^{t,\mu}) +  \nonumber \\
& & \;       \E \Big[  \partial_\mu w(s,\P_s^{t,\mu})(X_s^{t,\xi}).
b(s,X_s^{t,\xi},\tilde\alpha(s,X_s^{t,\xi},\P_s^{t,\mu}),Id\tilde\alpha(s,.,\P_s^{t,\mu})\star\P_s^{t,\mu}) \nonumber \\
& & \;\;\;  + \; \frac{1}{2}  {\rm tr}\big(\partial_x\partial_\mu w(s,\P_s^{t,\mu})(X_s^{t,\xi}) \sigma_s\sigma_s\trans
(s,X_s^{t,\xi},\tilde\alpha(s,X_s^{t,\xi},\P_s^{t,\mu}),Id\tilde\alpha(s,.,\P_s^{t,\mu})\star\P_s^{t,\mu}) \big) \Big] ds \nonumber \\
&=&   w(t,\mu) \; + \;  \int_t^T  \Dt{w}(s,\P_s^{t,\mu})  \; + \;   <\Lc_s^{\tilde\alpha(s,.,\P_s^{t,\mu})} w(s,\P_s^{t,\mu}) ,\P_s^{t,\mu} > ds,  \label{Itow}
\enq
where we used in the second equality the fact that $\P_s^{t,\mu}$ is the distribution of $X_s^{t,\xi}$ for $s$ $\in$ $[t,T]$.   Since $x$ $\mapsto$ $\tilde\alpha(s,.,\P_s^{t,\mu})$ $\in$ $L(\R^d;A)$ for $s$ $\in$ $[t,T]$, we deduce from the Bellman equation satisfied by $w$ and \reff{Itow}  that
\beqs
\hat g(\P_T^{t,\mu})  & \geq & w(t,\mu) - \int_t^T  \hat f(s,\P_s^{t,\mu},\tilde\alpha(s,.,\P_s^{t,\mu}))ds.
\enqs
Since $\alpha$ is arbitrary in $\Ac$, this shows that $w(t,\mu)$ $\leq$ $v(t,\mu)$.

In the final step, let us apply the same It\^o's argument \reff{Itow} with the feedback control $\alpha^*$ $\in$ $\Ac$ associated with the fonction
$\tilde\alpha^*$ $\in$ $Lip([0,T]\times\R^d\times\Pc_{_2}(\R^d);A)$.  Since $\tilde\alpha$ attains the infimum in \reff{HJB}, we thus get
\beqs
\hat g(\P_T^{t,\mu})  & =  & w(t,\mu) - \int_t^T  \hat f(s,\P_s^{t,\mu},\tilde\alpha^*(s,.,\P_s^{t,\mu}))ds,
\enqs
which shows that $w(t,\mu)$ $=$ $J(t,\mu,\alpha^*)$ ($\geq$ $v(t,\mu)$),
and therefore gives the required result: $v(t,\mu)$ $=$ $w(t,\mu)$ $=$ $J(t,\mu,\alpha^*)$.
\ep

\vspace{3mm}

We  shall apply the verification theorem in the next section, where we can derive  explicit (smooth) solutions to the Bellman equation \reff{HJB} in some class of  examples, but first discuss below the case when there are no mean-field interaction, and the structure of the optimal control (when it exists).

\begin{Remark} \label{remnomean}
{\rm (No mean-field interaction)

\noindent 
We consider the classical case of stochastic control  where there is no mean-field interaction in the dynamics of the state process, i.e. $b(t,x,a)$ and
$\sigma(t,x,a)$ do not depend on $\lambda$, as well as in the cost functions  $f(t,x,a)$ and $g(x)$. In this special case, let us show how the verification Theorem \ref{theoverif} is reduced to the  
classical verification result for smooth functions on $[0,T]\times\R^d$, see e.g. \cite{FleSon06} or \cite{Pha09}.

Suppose that there exists a function $u$ in $\Cc^{1,2}([0,T]\times\R^d)$ solution to the standard HJB equation
\begin{equation}\label{HJBstandard}
\left\{
\begin{array}{rcl}
\partial_t u + \Inf_{a\in A} \big[ f(t,x,a) + L_t^a u(t,x)  \big] &=& 0,  \; \mbox{ on } \; [0,T)\times\R^d,  \\
u(T,.) &=& g \;\; \mbox{ on } \;  \R^d.
\end{array}
\right.
\end{equation}
where $L_t^a$ is the second-order differential operator
\beqs
L_t^a u(t,x) &=& \partial_x u(t,x). b(t,x,a) + \frac{1}{2}{\rm tr}\big(\partial_{xx}^2  u(t,x) \sigma\sigma\trans(t,x,a) \big),
\enqs
and  that for all $(t,x)$ $\in$ $[0,T)\times\R^d$, there exists $\hat a(t,x)$ attaining the argmin in \reff{HJBstandard}, s.t. the map
$x$ $\mapsto$ $\hat a(t,x)$ is Lipschitz on $\R^d$.

Let us then consider the function defined on $[0,T]\times\Pc_{_2}(\R^d)$ by
\beqs
w(t,\mu) &=& <u(t,.),\mu> \; = \; \int_{\R^d} u(t,x) \mu(dx).
\enqs
The lifted function of $w$ is thus equal to $\Wc(t,X)$ $=$ $\E[u(t,X)]$ with Fr\'echet derivative (with respect to $X$ $\in$ $L^2(\Fc_0,\P))$:
$[D\Wc](t,X)(Y)$ $=$ $\E[\partial_x u(t,X).Y]$. Assuming that the time derivative of $u$ w.r.t. $t$ satisfies a quadratic growth condition in $x$, 
the first derivative  of $u$ w.r.t. $x$ satisfies a linear growth condition, and the second
derivative of $u$ w.r.t. $x$  is bounded, this shows that $w$ lies in $\Cc_b^{1,2}([0,T]\times\Pc_{_2}(\R^d))$ with
\beqs
\partial_t w(t,\mu) \; = \; <\partial_t u(t,.),\mu>, & &   \partial_\mu w(t,\mu) \; = \; \partial_x u(t,.), \;\;\;
\partial_x\partial_\mu v(t,\mu) \; = \;  \partial^2_{xx} u(t,.). 
\enqs
Recalling the definition \reff{Lcalpha} of $\Lc_t^{\tilde\alpha}w(t,\mu)$, we then get for any fixed $(t,\mu)$ $\in$ $[0,T)\times\Pc_{_2}(\R^d)$:
\beq
& & \partial_t w(t,\mu)  \; + \;   \Inf_{\tilde\alpha\in L(\R^d;A)} \Big[ \hat f(t,\mu,\tilde\alpha)  \; + \;   <\Lc_t^{\tilde\alpha} w(t,\mu) , \mu> \Big] \nonumber \\
&=& \Inf_{\tilde\alpha\in L(\R^d;A)}  \int_{\R^d} \big[ \partial_t u(t,x)  + f(t,x,\tilde\alpha(x)) + L_t^{\tilde\alpha(x)} u(t,x) \big] \mu(dx) \nonumber \\
&=& \int_{\R^d} \inf_{a\in A} \big[ \partial_t u(t,x)  + f(t,x,a) + L_t^{a} u(t,x) \big] \mu(dx).  \label{egalHJB}
\enq
Indeed, the inequality $\geq$ in \reff{egalHJB} is clear since $\tilde\alpha(x)$ lies in $A$  for all $x$ $\in$ $\R^d$, and $\tilde\alpha$ $\in$ $L(\R^d;A)$.
Conversely, by taking $\hat a(t,x)$ which attains the infimum  in \reff{HJBstandard}, and since the map $x$ $\in$ $\R^d$
$\mapsto$ $\hat a(t,x)$ is Lipschitz, we then have
\beqs
&& \int_{\R^d} \inf_{a\in A} \big[ \partial_t u(t,x)  + f(t,x,a) + L_t^{a} u(t,x) \big] \mu(dx) \\
&=&  \int_{\R^d}  \big[ \partial_t u(t,x)  + f(t,x,\hat a(t,x)) + L_t^{\hat a(t,x)} u(t,x) \big] \mu(dx) \\
& \geq & \Inf_{\tilde\alpha\in L(\R^d;A)}  \int_{\R^d} \big[ \partial_t u(t,x)  + f(t,x,\tilde\alpha(x)) + L_t^{\tilde\alpha(x)} u(t,x) \big] \mu(dx),
\enqs
which thus shows the equality \reff{egalHJB}. Since $u$ is solution to \reff{HJBstandard}, this proves that $w$ is solution to the Bellman equation
\reff{HJB},  $\tilde\alpha^*(t,x)$ $=$ $\hat a(t,x)$ is an optimal feedback control, and therefore, the value function is equal to
$v(t,\mu)$ $=$ $<u(t,.),\mu>$.
\ep
}
\end{Remark}

\begin{Remark}
{\rm (Form of the optimal control) 

\noindent Consider the case where the coefficients of the McKean-Vlasov SDE and of the running costs do not depend upon the law of the control, hence  in the form: $b(t,X_t,\alpha_t,\P_{_{X_t}})$, $\sigma(t,X_t,\alpha_t,\P_{_{X_t}})$, $f(t,X_t,\alpha_t,\P_{_{X_t}})$, and denote by 
\beqs
\H(t,x,a,\mu,q,M) &=& f(t,x,a,\mu) + q.b(t,x,a,\mu)  + \frac{1}{2}{\rm tr}\big(M \sigma\sigma\trans(t,x,a,\mu) \big)
\enqs
for $(t,x,a,\mu,q,M)$ $\in$ $[0,T]\times\R^d\times A\times\Pc_{_2}(\R^d)\times\R^d\times\S^{d}$, the Hamiltonian function related to the Bellman equation \reff{HJB} rewritten as:
\beq \label{HJBinfi}
\partial_t w(t,\mu) \; + \;   
\Inf_{\tilde\alpha\in L(\R^d;A)} \int_{\R^d} \H\big(t,x,\tilde\alpha(x),\mu,\partial_\mu w(\mu)(x),\partial_x\partial_\mu w(\mu)(x) \big) \mu(dx) &=& 0, 
\enq  
for $(t,\mu)$ $\in$ $[0,T)\times\Pc_{_2}(\R^d)$. Under suitable convexity conditions on the function $a$ $\in$ $A$ $\mapsto$ $\H(t,x,a,\mu,q,M)$, there exists a minimizer, say $\hat a(t,x,\mu,q,M)$, to $\inf_{a\in A} \H(t,x,a,\mu,q,M)$.  Then, an optimal control $\tilde\alpha^*$ in the statement of the verification theorem \ref{theoverif}, obtained from the minimization of the (infinite dimensional) Hamiltonian in \reff{HJBinfi}, is written merely as 
$\tilde\alpha^*(t,x,\mu)$ $=$ $\hat a(t,x,\mu,\partial_\mu w(\mu)(x),\partial_x\partial_\mu w(\mu)(x))$, which extends the form discuss in Remark \ref{remnomean}, and says that it depends locally upon the derivatives of the value function.  In the more general case when the coefficients depend upon the law of the control, we shall see how one can derive the form of the optimal control for the linear-quadratic problem. 
}
\ep
\end{Remark}

%However, in general, there is no smooth solutions to the HJB equation, and in the spirit of HJB equation for standard stochastic control,
%we shall introduce in the next paragraph  a notion of viscosity solutions for the Bellman equation \reff{HJB} in $[0,T]\times\Pc_{_2}(\R^d)$.

\section{Application: linear-quadratic McKean-Vlasov control pro\-blem} \label{secappli}

\setcounter{equation}{0} \setcounter{Assumption}{0}
\setcounter{Theorem}{0} \setcounter{Proposition}{0}
\setcounter{Corollary}{0} \setcounter{Lemma}{0}
\setcounter{Definition}{0} \setcounter{Remark}{0}

We consider a multivariate linear McKean-Vlasov controlled dynamics with coefficients given by
\begin{equation} \label{bsigLQ}
%\left\{
\begin{array}{ccl}
b(t,x,\mu,a,\lambda) &=&  b_0(t) + B(t) x + \bar B(t) \bar\mu + C(t) a + \bar C(t) \bar\lambda, \\
\sigma(t,x,\mu,a,\lambda) &=& \sigma_0(t) +  D(t) x + \bar D(t) \bar\mu + F(t) a + \bar F(t) \bar\lambda,
\end{array}
%\right.
\end{equation}
for $(t,x,\mu,a,\lambda)$ $\in$ $[0,T]\times\R^d\times\Pc_{_2}(\R^d)\times\R^m\times\Pc_{_2}(\R^m)$, where we set
\beqs
\bar\mu \; := \; \int_{\R^d} x \mu(dx), & & \bar\lambda \; := \; \int_{\R^m} a \lambda(da). 
\enqs
%\beqs
%dX_t &=& \big(B(t)X_s + \bar B(t) \E[X_t] + C(t) \alpha_t + \bar C(t) \E[\alpha_t] \big) dt  \\
%& & \;\;\; + \; \big(  D(t) X_t + \bar D(t) \E[X_t] + F(t) \alpha_t  + \bar F(t) \E[\alpha_t] \big) dB_t,
%\enqs
Here  $B$, $\bar B$,  $D$, $\bar D$ are deterministic continuous functions valued in $\R^{d\times d}$, and $C$, $\bar C$, $F$, $\bar F$ are deterministic continuous functions valued in $\R^{d\times m}$, and $b_0$, $\sigma_0$ are deterministic continuous function valued in $\R^d$. 
The quadratic cost functions  are given by
\begin{equation} \label{fgLQ}
%\left\{
\begin{array}{rcl}
f(t,x,\mu,a,\lambda) &=& x\trans Q_2(t)x + \bar\mu\trans \bar Q_2(t) \bar \mu  + a\trans R_2(t) a + \bar\lambda\trans \bar R_2(t) \bar\lambda  +  
2 x\trans M_2(t)a \\
& &  \;  + \;  2 \bar\mu\trans\bar M_2(t)\bar\lambda +  q_1(t).x + \bar q_1(t).\bar\mu + r_1(t).a + \bar r_1(t).\bar\lambda, \\
g(x,\mu) &=&  x\trans P_2 x + \bar\mu\trans \bar P_2 \bar\mu + p_1.x + \bar p_1.\bar\mu,
\end{array}
%\right.
\end{equation}where $Q_2$, $\bar Q_2$  are deterministic continuous functions, $P_2$, $\bar P_2$ are constants valued in $\R^{d\times d}$, $R_2$, $\bar R_2$ are deterministic continuous  functions valued in $\R^{m\times m}$,  $M_2$, $\bar M_2$ are deterministic continuous functions valued in $\R^{d\times m}$,  
$q_1$, $\bar q_1$ are deterministic continuous functions, $p_1$, $\bar p_1$ are constants valued in $\R^d$, and $r_1$, $\bar r_1$ are deterministic continuous  functions valued in $\R^m$.  Since $f$ and $g$ are real-valued, 
we may assume w.l.o.g. that all the matrices $Q_2$, $\bar Q_2$,  $R_2$, $\bar R_2$, $P_2$, $\bar P_2$ are symmetric. We denote by  $\S_+^d$ the set of nonnegative symmetric matrices in $\S^d$, and by $\S_{>+}^d$ the subset of symmetric  positive definite matrices.  
This linear quadratic (LQ) framework is similar to the one in \cite{yon13}, and extends the one considered in \cite{benetal16} where there is no dependence on the law of the control, and the diffusion coefficient is deterministic.

The functions $\hat f$ and $\hat g$ defined in \reff{defhatfg} are then given by
\begin{equation} \label{hatfgLQ}
\left\{
\begin{array}{rcl}
\hat f(t,\mu,\tilde\alpha) &=& {\rm Var}(\mu)(Q_2(t)) + \bar\mu \trans(Q_2(t) + \bar Q_2(t))\bar\mu\\
 & &  \;\;+\;\;{\rm Var}(\tilde\alpha \star \mu)(R_2(t)) + \overline{\tilde\alpha\star\mu} \trans (R_2(t) +\bar R_2(t)) \overline{\tilde\alpha \star \mu}\\
 & & \;\; + \;\;2 \bar \mu \trans(M_2(t) + \bar M_2(t))\overline{\tilde\alpha \star \mu} \; + \;  2 \int_{\R^d}  (x-\bar \mu)\trans M_2(t) \tilde\alpha(x) \mu(dx) \\
 & & \;\; +\;\;\big(q_1(t) + \bar q_1(t)\big).\bar\mu + \big(r_1(t) + \bar r_1(t)\big).\overline{\tilde\alpha \star\mu} \\
 \hat g(\mu) &=& {\rm Var}(\mu)(P_2) + \bar\mu \trans(P_2 + \bar P_2)\bar\mu + (p_1 + \bar p_1).\bar\mu,
\end{array}
\right.
\end{equation}
for any $(t,\mu)$ $\in$ $[0,T)\times\Pc_{_2}(\R^d)$, $\tilde\alpha$ $\in$ $L(\R^d;A)$ (here with $A$ $=$ $\R^m$),  
where  we set for any  $\Lambda$ in $\S^{d}$ (resp. in $\S^{m}$), and 
$\mu$ $\in$ $\Pc_{_2}(\R^d)$ (resp. $\Pc_{_2}(\R^m)$): 
\beqs
\bar\mu_{_2}(\Lambda)  \; := \;  \int x\trans \Lambda x \mu(dx),  & & {\rm Var}(\mu)(\Lambda) \; := \; \bar\mu_{_2}(\Lambda)  - \bar\mu\trans\Lambda\bar\mu.
\enqs

We look for a value function solution to the Bellman equation \reff{HJB} in the form
\beq \label{wquadra}
w(t,\mu) &=& {\rm Var}(\mu)(\Lambda(t)) + \bar\mu\trans\Gamma(t)\bar\mu + \gamma(t).\bar\mu + \chi(t),
\enq
for some  functions $\Lambda$, $\Gamma$ $\in$ $\Cc^1([0,T];\S^d)$, $\gamma$ $\in$ $\Cc^1([0,T];\R^d)$, and $\chi$ $\in$ $\Cc^1([0,T];\R)$.  The lifted function of $w$ in \reff{wquadra} is given by 
\beqs
\Wc(t,X) &=& \E[X\trans\Lambda(t)X] +  \E[X]\trans(\Gamma(t)-\Lambda(t))\E[X]  + \gamma(t).\E[X] + \chi(t),
\enqs
for $X$ $\in$ $L^2(\Fc_0;\R^d)$.  By computing for all $Y$ $\in$ $L^2(\Fc_0;\R^d)$ the difference
\beqs
\Wc(t, X+Y)-\Wc(t, X) &=&
%\E[(X+Y)\trans \Lambda(t) (X+Y)] -\E[X \trans \Lambda(t) X]\\
% & & \;\;- \; \E[X+Y]\trans(\Lambda(t)  - \Gamma(t))\E[X+Y] + \E[X]\trans(\Lambda(t) - \Gamma(t))\E[X]\\
% & &\;\; +\; \gamma(t). (\E[X+Y]-\E[X])\\
%&=& \E[Y \trans \Lambda(t) X + X \trans \Lambda(t) Y] -  \E[Y] \trans (\Lambda(t) - \Gamma(t)) \E[X] -  \E[X]\trans (\Lambda(t) - \Gamma(t))\E[Y] \\
%& &\;\; +\; \gamma(t) . \E[Y] + o(\vert Y \vert_{L^2}) \\
%&=&2 \E[X \trans \Lambda(t) Y] (+?) \red{-}   2\E[X] \trans(\Lambda(t)- \Gamma(t)) \E[Y] + \gamma(t) \E[Y] + o(\vert Y \vert_{L^2})\\
%&=& 
\E\Big[\big(2X \trans \Lambda(t) +  2 \E[X]\trans(\Gamma(t) - \Lambda(t)) + \gamma(t)\big).Y\Big] \; + \; o(  \|Y \|_{_{L^2}}),
\enqs
we see that $\Wc$ is Fr\'echet differentiable (w.r.t. $X$) with $[D\Wc](t,X)(Y)$ $=$ $\E\big[\big(2X \trans \Lambda(t) +  2 \E[X]\trans(\Gamma(t) - \Lambda(t)) + \gamma(t)\big).Y\big]$.  
This shows that $w$ lies in $\Cc_b^{1,2}([0,T]\times\Pc_{_2}(\R^d))$ with
\beqs
\partial_t w(t, \mu)  &=& {\rm Var}(\mu) (\Lambda'(t)) +\bar\mu\trans\Gamma'(t)\bar\mu + \gamma'(t)\bar\mu + \chi'(t), \\
\partial_\mu w(t, \mu)(x)&=& 2 x \trans \Lambda(t) + 2 \bar \mu \trans(\Gamma(t)-\Lambda(t)) + \gamma(t), \\
\partial_x\partial_\mu w(t, \mu)(x)&=& 2 \Lambda(t).
\enqs
Together with the quadratic expression \reff{hatfgLQ} of $\hat f$, $\hat g$, we then see that $w$ satisfies the Bellman equation \reff{HJB} iff
\beq
& & {\rm Var}(\mu)(\Lambda(T)) + \bar\mu\trans\Gamma(T)\bar\mu + \gamma(T).\bar\mu + \chi(T) \nonumber \\
&=&  {\rm Var}(\mu)(P_2) + \bar\mu \trans(P_2 + \bar P_2)\bar\mu + (p_1 + \bar p_1).\bar\mu, \label{HJBLQT} 
\enq
holds for all  $\mu$ $\in$ $\Pc_{_2}(\R^d)$, and 
\beq
& & {\rm Var} (\mu)\big(\Lambda'(t)+ Q_2(t)+ D(t)\trans \Lambda(t)D(t) + \Lambda(t)B(t)+ B(t)\trans \Lambda(t)\big) 
+  \inf_{\tilde\alpha \in L(\R^d, A)} G_t^{\mu}(\tilde\alpha) \nonumber \\
& & \;\;+\; \bar\mu\trans\Big(\Gamma'(t) +Q_2(t) +\bar Q_2(t)+ (D(t) + \bar D(t)) \trans \Lambda(t)(D(t) +\bar D(t)) \nonumber \\
& & \hspace{2cm}  + \;  \Gamma(t)(B(t)+ \bar B(t)) + (B(t) \;+\; \bar B(t))\trans \Gamma(t) \Big) \bar\mu \nonumber \\
& & \; \; +  \;  \big( q_1(t) + \bar q_1(t)+ \gamma(t)(B(t)+\bar B(t))    + 2 \sigma_0\trans \Lambda(t) (D(t) + \bar D(t))  +  2 b_0(t)\trans \Gamma(t) \big)\bar\mu  \nonumber \\
& &  \; +\; \chi'(t)   \; + \;  \gamma(t).b_0(t) + \sigma_0(t)\trans \Lambda(t)b_0(t) \nonumber \\
& =&  0, \label{HJBLQt}
%& &+(\gamma'(t)+ q_1(t) +\bar q_1(t)+ \gamma(t)(B(t)+ \bar B(t)))\bar\mu +\chi'(t) =0
\enq
holds for all  $t$ $\in$ $[0,T)$, $\mu$ $\in$ $\Pc_{_2}(\R^d)$, where the function $G_t^\mu$ $:$ $L_\mu^2(A)$ $\supset$ $L(\R^d;A)$  $\rightarrow$ $\R$ is defined by  
\beq 
G_t^\mu(\tilde\alpha)&=& {\rm Var}(\tilde\alpha \star \mu) (U_t) \;+\; \overline{\tilde\alpha \star \mu}\trans V_t   \overline{\tilde\alpha \star \mu} 
\; +\; 2  \int_{\R^d} (x - \bar \mu)\trans S_t  \tilde\alpha(x)\mu(dx)  \nonumber \\
 & & \;+\; 2\bar\mu\trans Z_t \overline{\tilde\alpha\star \mu}   \;+\;  Y_t . \overline{\tilde\alpha\star \mu}, \label{defGk} 
 \enq
 and we set  $U_t$ $=$ $U(t,\Lambda(t))$, $V_t$ $=$ $V(t,\Lambda(t))$, $S_t$ $=$ $S(t,\Lambda(t))$, $Z_t$ $=$ $Z(t,\Lambda(t),\Gamma(t))$, 
 $Y_t$ $=$ $Y(t,\Gamma(t),\gamma(t))$ with
\begin{equation} \label{defUVSZ}
\left\{
\begin{array}{rcl}
 U(t,\Lambda(t)) &=&F(t)\trans \Lambda(t) F(t)+ R_2(t),\\
 V(t,\Lambda(t)) &=&(F(t) + \bar F(t))\trans \Lambda(t) (F(t)+\bar F(t))+ R_2(t) + \bar R_2(t),\\
 S(t,\Lambda(t)) &=&D(t)\trans \Lambda(t) F(t) + \Lambda(t) C(t) + M_2(t),\\
 Z(t,\Lambda(t),\Gamma(t)) &=& (D(t) + \bar D(t)) \trans \Lambda(t)(F(t) + \bar F(t)) + \Gamma(t)(C(t) + \bar C(t)) + M_2(t) + \bar M_2(t) \\
 Y(t,\Gamma(t),\gamma(t)) &=& \big(C(t)+ \bar C(t)\big)\trans\gamma(t) + r_1(t) + \bar r_1(t)+ 2 \big(F(t)+ \bar F(t)\big)\trans \Lambda(t)\sigma_0(t). 
\end{array}
\right. 
\end{equation}
%Here, $L_\mu^2(\R^d;A)$ $\supset$ $L(\R^d;A)$ is the Hilbert space of measurable functions on  $\R^d$ valued in $A$ $=$ $\R^m$ 
 %and square integrable w.r.t. $\mu$ $\in$ $\Pc_{_2}(\R^d)$. 
 %We now search for the infimum of the function $G_t^{\mu}$, and shall make the following assumptions on the  symmetric matrices of the quadratic cost functional:
 %\begin{equation} \label{condmatrice}
%\left\{
%\begin{array}{ccl}
 %P_2  \; \geq \; 0, \; P_2 +\bar P_2  \; \geq \; 0, & & Q_2(t) \; \geq \; 0, \; Q_2(t) + \bar Q_2(t) \; \geq \; 0,  \\
 %& & R_2(t) \; >  \; 0, \; R_2(t) + \bar R_2(t) \; >  \;0
%\end{array}
%\right.
%\end{equation}
We now search for the infimum of the function $G_t^{\mu}$. After some  straightforward calculation,  we derive the Gateaux derivative of $G_t^\mu$ at $\tilde\alpha$ in the direction $\beta$ $\in$ $L_\mu^2(A)$, which  is given by:
\beqs
DG_t^{\mu}(\tilde\alpha, \beta)&: =& \lim_{\eps\rightarrow 0}  \frac{G_t^\mu(\tilde\alpha + \eps\beta) - G_t^\mu(\tilde\alpha)}{\eps}
\; = \; \int_{\R^d} \dot g_t^\mu(x,\tilde\alpha).\beta(x) \mu(dx)
\enqs
% & & \;\;+\; \big(r_1(t) + \bar r_1(t)+ \gamma(t)(C(t) +\bar C(t))\big) \tilde\alpha\star\mu
with
\beqs
\dot g_t^\mu(x, \tilde\alpha) &=& 2 U_t \tilde\alpha + 2 (V_t -U_t) \overline{\tilde\alpha\star\mu} 
+  2S_t\trans(x-\bar\mu) + 2 Z_t\trans\bar\mu  + Y_t.
\enqs
%We shall see later than under \reff{condmatrice}, the symmetric matrices $U_t$ and $V_t$ in \reff{defUVSZ} are positive, hence invertible. 
Suppose that the symmetric matrices $U_t$ and $V_t$ in \reff{defUVSZ} are positive, hence invertible (this will be discussed later on).  Then, the function $G_t^\mu$ is convex and coercive on the Hilbert space $L_\mu^2(A)$, and attains its infimum 
at some $\tilde\alpha$ $=$  $\tilde\alpha^*(t,.,\mu)$ s.t.  $DG_t^\mu(\tilde\alpha;.)$ vanishes, i.e.  $\dot g_t^\mu(x,\tilde\alpha^*(t,.,\mu))$ $=$ $0$ for all $x$ $\in$ $\R^d$, which gives:
\beq
\tilde\alpha^*(t,x,\mu)  &=&   - U_t^{-1} S_t\trans(x-\bar\mu)\;  -\;   V_t^{-1} Z_t \trans \bar\mu \;-\; 
\frac{1}{2} V_t^{-1}Y_t.   \label{optialphaLQ}
\enq
It is clear  that  $\tilde\alpha^*(t,.,\mu)$ lies in $L(\R^d;A)$,  and so  after some straightforward caculation:
\beqs
\inf_{\tilde\alpha \in L(R^d, A)} G_t^{\mu}(\tilde\alpha) \; = \; G_t^\mu(\tilde\alpha^*(t,.,\mu))&= & - {\rm Var}(\mu)\big(  S_t  U_t^{-1} S_t \trans \big)  
 \;  - \;  \bar\mu\trans \big( Z_t V_t^{-1}Z_t\trans\big) \bar\mu  \\
& & \;\; - \; Y_t\trans V_t^{-1}Z_t \trans \bar\mu  \; -\; \frac{1}{4}  Y_t\trans V_t^{-1} Y_t. 
\enqs
Plugging the above expression in \reff{HJBLQt}, we observe that  the relation  \reff{HJBLQT}-\reff{HJBLQt}, hence the Bellman equation,  
is satisfied by identifying the terms in ${\rm Var}(\mu)(.)$, $\bar\mu\trans(.)\bar\mu$,  $\bar\mu$,  which leads to the system of ordinary differential 
equations (ODEs) for  $(\Lambda,\Gamma,\gamma,\chi)$:   
\begin{equation} \label{Riccatilambda}
\left\{
\begin{array}{rcl}
\Lambda'(t)+ Q_2(t)+ D(t)\trans \Lambda(t)D(t) +\Lambda(t)B(t) + B(t) \trans \Lambda(t) && \\
 \;\;\; -S(t,\Lambda(t))U(t,\Lambda(t))^{-1} S(t,\Lambda(t))\trans &= & 0,\\
\Lambda(T) &=& P_2,
\end{array}
\right.
\end{equation}
\begin{equation}\label{Riccatigamma}
\left\{
\begin{array}{rcl}
\Gamma'(t) +Q_2(t) +\bar Q_2(t)+ (D(t) + \bar D(t)) \trans \Lambda(t)(D(t) +\bar D(t))+ \Gamma(t)(B(t)+ \bar B(t))  \\
\;\;+\; (B(t) + \bar B(t))\trans \Gamma(t)- Z(t,\Lambda(t),\Gamma(t))  V(t,\Lambda(t))^{-1} Z(t,\Lambda(t),\Gamma(t))\trans &=& 0,\\
\Gamma(T) \; = \;  P_2 + \bar P_2, & & 
\end{array}
\right.
\end{equation}
\begin{equation} \label{lingamma}
\left\{
\begin{array}{rcl}
\gamma'(t) +  \big(B(t) + \bar B(t))\trans \gamma(t)  -  Z(t,\Lambda(t),\Gamma(t)) V(t,\Lambda(t))^{-1} Y(t,\Gamma(t),\gamma(t))  \\
\;\;\; + \;  q_1(t) + \bar q_1(t) + 2 \big( D(t) + \bar D(t) \big)\trans \Lambda(t) \sigma_0(t) + 2 \Gamma(t) b_0(t)  &=& 0,\\
\gamma(T) \; = \;  p_1 + \bar p_1 & & 
\end{array}
\right.
\end{equation}
\begin{equation} \label{linchi}
\left \{
\begin{array}{rcl}
\chi'(t) -\; \frac{1}{4} Y(t,\Gamma(t),\gamma(t))\trans V(t,\Lambda(t))^{-1}Y(t,\Gamma(t),\gamma(t)) \\
 \;\;\; +  \; \gamma(t).b_0(t)  + \sigma_0(t)\trans\Lambda(t)\sigma_0(t)  &=& 0, \\
\chi(T) \; = \;  0. & & 
\end{array}
\right.
\end{equation}

Therefore, the resolution of the Bellman equation in the LQ framework is reduced to the resolution of the Riccati equations  \reff{Riccatilambda}  and  \reff{Riccatigamma} for $\Lambda$ and $\Gamma$, and then 
given $(\Lambda,\Gamma)$, to the resolution  of the  linear ODEs \reff{lingamma} and \reff{linchi} for $\gamma$ and $\chi$.  Suppose that there exists a solution $(\Lambda,\Gamma)$ $\in$ 
$\Cc^1([0,T];\S^d)\times \Cc^1([0,T];\S^d)$  to \reff{Riccatilambda}-\reff{Riccatigamma}  
s.t. $(U_t,V_t)$ in \reff{defUVSZ} lies in $\S^m_{>+}\times\S^m_{>+}$  for all $t$ $\in$ $[0,T]$ (see Remark \ref{remriccati}). Then,  the above calculations 
are justified a posteriori, and by noting also that  the mapping $(t,x,\mu)$ $\mapsto$ $\tilde\alpha^*(t,x,\mu)$ $\in$  $Lip([0,T]\times\R^d\times\Pc_{_2}(\R^d);A)$, 
we deduce by the verification theorem that the value function $v$ is equal to $w$ in \reff{wquadra} with $(\Lambda,\Gamma,\gamma,\chi)$ solution to 
\reff{Riccatilambda}-\reff{Riccatigamma}-\reff{lingamma}-\reff{linchi}.  Moreover, the optimal control is given in feedback form from \reff{optialphaLQ} by 
\beq \label{opticontrolLQ}
\alpha_t^* &=& \tilde\alpha^*(t,X_t^*,\P_{_{X_t^*}})   =  
 - U_t^{-1} S_t\trans(X_t^* - \E[X_t^*])\;  -\;   V_t^{-1} Z_t \trans \E[X_t^*] \;-\;  \frac{1}{2} V_t^{-1}Y_t, 
\enq
where $X^*$ is the state process controlled by $\alpha^*$.  
%EXPLICIT THE FORM OF THE CONTROL AS IN THE DISCRETE TIME CASE BY COMPUTING $\E[X_t^*]$ 

\vspace{2mm}

\begin{Remark} \label{remriccati}
{\rm In the case where $M_2$ $=$ $\bar M_2$ $=$ $0$  (i.e. no crossing term between the state and the control in the quadratic cost function $f$), it is shown in  Proposition 3.1 and 3.2 in  
\cite{yon13} that  under the condition 
\begin{equation} \label{condmatrice}
%\left\{
\begin{array}{ccl}
 P_2  \; \geq \; 0, \; P_2 +\bar P_2  \; \geq \; 0, & & Q_2(t) \; \geq \; 0, \; Q_2(t) + \bar Q_2(t) \; \geq \; 0,   \\
 & & R_2(t) \; \geq  \delta I_m, \; R_2(t) + \bar R_2(t) \; \geq   \; \delta I_m 
\end{array}
%\right.
\end{equation}
for some $\delta$ $>$ $0$,   the Riccati equations \reff{Riccatilambda}-\reff{Riccatigamma}  admit unique solutions  $(\Lambda,\Gamma)$ $\in$  
$\Cc^1([0,T];\S_+^d)$ $\times$ $\Cc^1([0,T];\S_+^d)$,  and then $U_t,V_t$ in \reff{defUVSZ} are symmetric positive definite matrices, i.e. 
lie in $\S_{>+}^m$ for all $t$ $\in$ $[0,T]$. 
In this case, we retrieve the expressions \reff{opticontrolLQ} of the optimal control in feedback form obtained in \cite{yon13}. 

We shall see in the next two paragraphs, some other examples arising from finance  with explicit solutions where condition \reff{condmatrice} is not satisfied.  
}
\ep
\end{Remark}

%Complete and compare with the results obtained in \cite{benetal13}, \cite{yon13} or \cite{cardel13}.

\vspace{3mm}

\subsection{Mean-variance portfolio selection} \label{subMV}

The  mean-variance problem consists in minimizing a cost functional  of the form:
 \beqs
 J(\alpha) &=& \frac{\eta}{2} {\rm Var}(X_T) -  \E[X_T]  \\
 &=& \E\big[ \frac{\eta}{2} \big(X_T\big)^2 -  X_T \big] \; - \;  \frac{\eta}{2} \Big( \E[X_T] \Big)^2
 \enqs
 for some $\eta$ $>$ $0$, with a dynamics for the wealth process $X$ $=$ $(X^\alpha)$  controlled by the amount $\alpha_t$ valued in $A$ $=$ $\R$
 invested in one risky stock at time $t$:
 \beqs
 dX_t &=& r(t) X_t dt + \alpha_t \big( \rho(t) dt + \vartheta(t) dB_t), \;\;\; X_0 \; = x_0 \; \in \;  \R, 
 \enqs
 where  $r$ is the interest rate, $\rho$  and $\vartheta$ $>$ $0$  are the excess rate of return (w.r.t. the interest rate) and volatility of the stock price, and these deterministic functions are assumed to be continuous.  
This model fits into the LQ framework \reff{bsigLQ}-\reff{fgLQ} of the McKean-Vlasov problem, with a linear controlled dynamics that does not have  mean-field  interaction: 
\beqs
b_0 = 0, \; B(t)=r(t), \;  \bar B =0, \; C(t)=\rho(t), \; \bar C =0,\\
\sigma_0 = D =  \bar D =0, \; F(t)=\vartheta(t), \; \bar F =0,\\
Q_2 =  \bar Q_2 =   M_2 = \bar M_2 =  R_2 =  \bar R_2 =0,\\
q_1 = \bar q_1 = r_1 = \bar r_1 =  0, \; P_2 = \frac{\eta}{2}, \; \bar P_2 =-\frac{\eta}{2}, \; p_1 = 0, \; \bar p_1 =-1.
\enqs
The Riccati system \reff{Riccatilambda}-\reff{Riccatigamma}-\reff{lingamma}-\reff{linchi} 
for $(\Lambda(t), \Gamma(t), \gamma(t), \chi(t))$ is written in this case as
\begin{equation}
\left\{
\begin{array}{rclrccc}
\Lambda'(t)-(\frac{\rho^2(t)}{\vartheta^2(t)}-2 r(t))\Lambda(t)&=&0, & &\Lambda(T) &=& \frac{\eta}{2},\\
\Gamma'(t) - \frac{\rho^2(t)\Gamma^2(t)}{\vartheta^2(t) \Lambda(t)}  + 2 r(t)\Gamma(t)  &=&0, & & \Gamma(T) &=& 0,\\
\gamma'(t) + r(t)\gamma(t)-\gamma(t)\frac{\rho^2(t)\Gamma(t)}{\vartheta^2(t)\Lambda(t)}&=&0, & & \gamma(T) &= & -1,\\
\chi'(t)-\frac{\rho^2(t)\gamma^2(t)}{4\vartheta^2(t)\Lambda(t)}&=&0, & &\chi(T)& =& 0,
\end{array}
\right.
\end{equation}
whose explicit solution is given by
\begin{equation}
\left\{
\begin{array}{ccl}
\Lambda(t)&=& \frac{\eta}{2}\exp\big(\int_t^T 2r(s)-\frac{\rho^2(s)}{\vartheta^2(s)} \; ds\big),\\
\Gamma(t)&=&0,\\
\gamma(t)&=&-\exp\big(\int_t^T r(s)ds\big)\\
\chi(t)&=&- \frac{1}{2\eta}\Big[ \exp\big(\int_t^T \frac{\rho^2(s)}{\vartheta^2(s)}ds\big)- 1 \Big].
\end{array}
\right.
\end{equation}
Although the condition \reff{condmatrice} is not satisfied, we see  that  $(U_t,V_t)$ in \reff{defUVSZ}, which are here explicitly given by 
$U_t$ $=$  $V_t$ $=$ $\vartheta(t)^2\Lambda(t)$, are positive, and this validates  our calculations for the verification theorem. Notice also that the functions $(Z_t,Y_t)$ in \reff{defUVSZ} are explicitly given by 
$Z_t$ $=$ $0$, $Y_t$ $=$ $\rho(t)\gamma(t)$. Therefore, the optimal  control is given in  feedback form from  \reff{opticontrolLQ} by 
\beq 
\alpha^*_t & = & \tilde\alpha^*(t,X_t^*, \P_{_{X_t^*}}) \nonumber \\
&=& -\frac{\rho(t)}{\vartheta^2(t)} (X_t^* -\E[X_t^*])   \; + \;    \frac{\rho(t)}{\eta \vartheta^2(t)} \exp\Big(\int_t^T \frac{\rho^2(s)}{\vartheta^2(s)}-r(s) \; ds\Big), \label{optimalcontrol}
\enq
where $X^*$ is the optimal wealth process with portfolio strategy $\alpha^*$, hence with mean process  governed by
\beqs \label{expectationofX}
d \E[X_t^*]&=&r(t) \E[X_t^*] dt  \;  + \;    \frac{\rho^2(t)}{\eta \vartheta^2(t)} \exp\Big( \int_t^T \frac{\rho^2(s)}{\vartheta^2(s)}-r(s) \; ds \Big) dt,
\enqs
and explicitly given by 
%The solution to \reff{expectationofX} is
\beqs \label{solution}
\E[X_t^*] &=& x_0  e^{\int_0^t r(s)ds}  \;+\; 
\frac{1}{\eta} \exp\Big(\int_t^T  \frac{\rho^2(s)}{\vartheta^2(s)}-r(s) \; ds\Big) \Big( \exp\big(\int_0^t \frac{\rho^2(s)}{\vartheta^2(s)}ds\big) -1 \Big). 
\enqs
Plugging into \reff{optimalcontrol}, we get the optimal control for the mean-variance portfolio problem 
\beqs
\alpha^*_t  &=& \frac{\rho(t)}{\vartheta^2(t)} \Big[ x_0 e^{\int_0^t r(s)ds} + \frac{1}{\eta} \exp\Big( \int_0^T  \frac{\rho^2(s)}{\vartheta^2(s)}  ds -  \int_t^T  r(s) ds\Big) - X_t^* \Big],
\enqs
and retrieve the closed-form expression of the optimal control  found in \cite{lizho00}, \cite{anddje10} or \cite{fisliv15} by different approaches.

\subsection{Inter-bank systemic risk model}

We consider a model of inter-bank borrowing and lending studied in \cite{caretal14} where the log-monetary reserve of each bank in the asymptotics when the number of banks tend to infinity, is governed by 
the McKean-Vlasov equation:
\beq \label{logX}
dX_t &=&  \big[ \kappa(\E[X_t] - X_t)  + \alpha_t ] dt + \sigma dB_t, \;\;\;\;\;  X_0 \; = \;  x_0  \in   \R.   
\enq 
Here, $\kappa$ $\geq$ $0$ is the rate of mean-reversion in the interaction from borrowing and lending between the banks, and $\sigma$ $>$ $0$ is the volatility coefficient of the bank reserve, assumed to be constant. 
Moreover, all  banks  can  control  their  rate of borrowing/lending to a central bank with the same policy $\alpha$ in order to minimize  a cost functional of the form 
\beqs
J(\alpha) &=&  \E \Big[ \int_0^T \Big( \frac{1}{2} \alpha_t^2 - q \alpha_t ( \E[X_t] - X_t) + \frac{\eta}{2} (\E[X_t] - X_t)^2 \Big) dt + \frac{c}{2}  (\E[X_T] - X_T)^2 \Big],   
\enqs
where $q$ $>$ $0$ is a positive parameter for the incentive to borrowing ($\alpha_t$ $>$ $0$) or lending ($\alpha_t$ $<$ $0$), and $\eta$ $>$ $0$, $c$ $>$ $0$ are positive parameters for penalizing departure from the average. 
This model fits into the LQ McKean-Vlasov  framework \reff{bsigLQ}-\reff{fgLQ}  with $d$ $=$ $m$ $=$ $1$ and
\beqs
b_0  =  0, \; B =  - \kappa, \; \bar B  =  \kappa, \; C  =  1, \; \bar C  =  0, \\
\sigma_0  =  \sigma, \; D  =  \bar D  =  F  =   \bar F =  0, \\
Q_2  =   \frac{\eta}{2}, \; \bar Q_2  =  - \frac{\eta}{2}, \;  R_2  =  \frac{1}{2}, \; \bar R_2  =  0, \; M_2  =  \frac{q}{2},\; \bar M_2  =  - \frac{q}{2}, \\
q_1 = \bar q_1 = r_1 = \bar r_1 =  0,  \;\;\; P_2 = \frac{c}{2}, \; \bar P_2 = - \frac{c}{2}, \; p_1 = \bar p_1 = 0. 
\enqs
The Riccati system \reff{Riccatilambda}-\reff{Riccatigamma}-\reff{lingamma}-\reff{linchi} 
for $(\Lambda(t), \Gamma(t), \gamma(t), \chi(t))$ is written in this case as
\begin{equation}
\left\{
\begin{array}{rclrccc}
\Lambda'(t) -2 (\kappa + q)\Lambda(t)  -  2 \Lambda^2 (t)  -\frac{1}{2}(q^2 -\eta) &=& 0, & &  \Lambda(T)&=& \frac{c}{2}, \\
\Gamma'(t) -2 \Gamma^2(t) &=& 0, & &  \Gamma(T)&=&0,\\
\gamma'(t)-2 \gamma(t) \Gamma(t) &=& 0, & & \gamma(T) &= &  0,\\
\chi'(t)+\sigma^2 \Lambda(t) -\frac{1}{2} \gamma^2(t) &=&0, & &  \chi(T) &= & 0,
\end{array}
\right.
\end{equation}
whose explicit solution is given by $\Gamma$ $=$ $\gamma$ $=$ $0$, and 
\beqs
\chi(t) &=& \sigma^2 \int_t^T \Lambda(s) ds,  \\
\Lambda(t) &=& \frac{1}{2} \frac{ (q-\eta^2)\big(  e^{(\delta^+-\delta^-)(T-t)}-1\big)  - c\big( \delta^+e^{(\delta^+-\delta^-)(T-t)} - \delta^- \big)}
{  \delta^- e^{(\delta^+-\delta^-)(T-t)} - \delta^+ \big)  - c e^{(\delta^+-\delta^-)(T-t)}-1}, 
\enqs 
where we set
\beqs
\delta^\pm &=& - (\kappa + q) \pm \sqrt{ (\kappa +q)^2 + (\eta-q^2)}. 
\enqs 
Moreover,  the functions $(U_t,V_t,Z_t,Y_t)$ in \reff{defUVSZ} are explicitly given by:   $U_t$ $=$ $V_t$ $=$ $\frac{1}{2}$ (hence $>$ $0$),  
$S_t$ $=$ $\Lambda(t)+\frac{q}{2}$, $Z_t$ $=$ $\Gamma(t)$ $=$ $0$, $Y_t$ $=$ $\gamma(t)$ $=$ $0$. Therefore,  the optimal  control is given in  feedback form from  \reff{opticontrolLQ} by 
\beq \label{controlsys}
\alpha_t^* \; = \; \tilde\alpha^*(t,X_t^*, \P_{_{X_t^*}})  &=& - (2\Lambda(t)+ q)(X_t^* -\E[X_t^*]),
\enq
where $X^*$ is the optimal  log-monetary reserve controlled by the rate  of borrowing/lending $\alpha^*$.  We then retrieve the expression found in \cite{caretal14} by sending the number of banks $N$ to infinity in their 
formula for the optimal control.  Actually, from \reff{logX}, we have $d\E[X_t^*]$ $=$ $\E[\alpha_t^*] dt$, while $\E[\alpha_t^*]$ $=$ $0$ from \reff{controlsys}.  We conclude that the optimal rate  of borro\-wing/lending is equal to
\beqs
\alpha_t^*  &=& - (2\Lambda(t)+ q)(X_t^* - x_0), \;\;\; 0 \leq t \leq T. 
\enqs

 \section{Viscosity solutions}

\setcounter{equation}{0} \setcounter{Assumption}{0}
\setcounter{Theorem}{0} \setcounter{Proposition}{0}
\setcounter{Corollary}{0} \setcounter{Lemma}{0}
\setcounter{Definition}{0} \setcounter{Remark}{0}

In general, there are no smooth solutions to the HJB equation, and in the spirit of HJB equation for standard stochastic control,
we shall introduce in this section a notion of viscosity solutions for the Bellman equation \reff{HJB} in the Wasserstein space of probability measures $\Pc_{_2}(\R^d)$.  
%There are several definitions of viscosity solutions for Hamilton Jacobi equations of first order in Wasserstein space and more generally in metric spaces, see e.g. \cite{ambetal05}, \cite{ganetal08}, 
%\cite{fenkat09} or \cite{ganswi15}.  
We adopt the approach in \cite{lio12}, and detailed in \cite{car12}, which consists, after the lifting identification between measures and random variables,  
in working in the Hilbert space $L^2(\Fc_0;\R^d)$  instead of working in  the Wasserstein space $\Pc_{_2}(\R^d)$, in order to use the various   tools developed for viscosity solutions in Hilbert spaces, in particular in our context,  
for second order Hamilton-Jacobi equations. 
%see \cite{lio88}, \cite{lio89a}, \cite{lio89b}, and the recent monograph \cite{fabgozswi15}.    

Let us rewrite the the Bellman equation \reff{HJB} in the ``Hamiltonian" form:
\begin{equation} \label{HJB2}
\left\{
\begin{array}{rcl}
- \Dt{v} + H(t,\mu,\partial_\mu v(t,\mu),\partial_x\partial_\mu v(t,\mu)) &=& 0    \;\;\; \mbox{ on } \; [0,T)\times\Pc_{_2}(\R^d), \\
v(T,.) &=& \hat g  \;\;\; \mbox{ on } \; \Pc_{_2}(\R^d) 
\end{array}
\right.
\end{equation}
 where $H$ is the function defined  by
\beq
 H(t,\mu,p,\Gamma) &=&-  \inf_{\tilde\alpha\in L(\R^d;A)} \Big[ < f(t,.,\mu,\tilde\alpha(.),Id\tilde\alpha\star\mu) + p(.).b(t,.,\mu,\tilde\alpha(.),Id\tilde\alpha\star\mu) \nonumber \\
& & \hspace{2cm}  + \; \frac{1}{2}{\rm tr}\big(\Gamma(.)\sigma\sigma\trans(t,.,\mu,\tilde\alpha(.),Id\tilde\alpha\star\mu) \big), \mu > \Big], \label{defH}
 \enq
 for  $(t,\mu)$ $\in$ $[0,T]\times\Pc_{_2}(\R^d)$, $(p,\Gamma)$ $\in$ $L^2_\mu(\R^d)\times L^\infty_\mu(\S^{d})$. 
 
We  then consider the ``lifted" Bellman equation in $ [0,T]\times L^2(\Fc_0;\R^d)$: 
\begin{equation} \label{HJBlift}
\left\{
\begin{array}{rcl}
- \Dt{V} + \Hc(t,\xi, DV(t,\xi),D^2 V(t,\xi)) &=& 0    \;\;\; \mbox{ on } \; [0,T)\times L^2(\Fc_0;\R^d), \\
V(T,\xi) &=& \hat G(\xi) \; := \; \E[g(\xi,\P_{_\xi})],   \;\;\;  \xi  \in   L^2(\Fc_0;\R^d), 
\end{array}
\right.
\end{equation}
where $\Hc$ $:$ $[0,T]\times L^2(\Fc_0;\R^d)\times L^2(\Fc_0;\R^d)\times S(L^2(\Fc_0;\R^d))$ $\rightarrow$ $\R$ is defined by 
\beq
\Hc(t,\xi,P,Q) &=& - \inf_{\tilde\alpha\in L(\R^d;A)} \Big\{ \E \Big[ f(t,\xi,\P_{_\xi},\tilde\alpha(\xi),Id\tilde\alpha\star\P_{_\xi})   +  P.b(t,\xi,\P_{_\xi},Id\tilde\alpha\star\P_{_\xi})  \label{defHc} \\
& & \hspace{0.9cm}  + \;  \frac{1}{2}  Q\big(\sigma(t,\xi,\P_{_\xi},\tilde\alpha(\xi),Id\tilde\alpha\star\P_{_\xi})N\big).\big(\sigma(t,\xi,\P_{_\xi},\tilde\alpha(\xi),Id\tilde\alpha\star\P_{_\xi}) N\big)  \Big] \Big\},  \nonumber 
\enq
with $N$ $\in$ $L^2(\Fc_0,\R^n)$ of zero mean, unit variance,  and independent of $\xi$.  
Observe that  when   $v$ and $V$ are smooth functions respectively in  $[0,T]\times\Pc_{_2}(\R^d)$ and $[0,T]\times L^2(\Fc_0;\R^d)$, linked by the lifting  relation $V(t,\xi)$ $=$ $v(t,\P_{_\xi})$, then from \reff{Uu1}-\reff{Uu2}, 
$v$ is solution to the Bellman equation \reff{HJB2} iff $V$ is solution to the Bellman equation \reff{HJBlift}.    
Let us mention that the lifted Bellman equation was also derived in \cite{benetal15} in the case where 
$\sigma$ $=$ $\sigma(x)$ is not controlled and does not depend on the distribution of the state process, and there is no dependence on the marginal distribution of the control process on the coefficients $b$ and $f$.

 \vspace{1mm}

 It is then natural to define viscosity solutions for the Bellman equation \reff{HJB2} (hence \reff{HJB})  from viscosity solutions to  \reff{HJBlift}. As usual, we say  that a function $u$ (resp. $U$) is locally bounded in 
 $[0,T]\times \Pc_{_2}(\R^d)$ (resp.  on  $[0,T]\times L^2(\Fc_0;\R^d)$) if it is bounded on bounded subsets of $[0,T]\times \Pc_{_2}(\R^d)$ (resp. of  $[0,T]\times L^2(\Fc_0;\R^d)$), and we denote by 
 $u^*$ (resp. $U^*$)  its upper   semicontinuous envelope, and by  $u_*$ (resp.  $U_*$) its lower semicontinuous envelope.  Similarly as in \cite{ganetal08}, we define the set 
 $\Cc_{\ell}^2([0,T]\times L^2(\Fc_0;\R^d))$   of test functions for the lifted Bellman equation, as  the set of real-valued continuous functions $\Phi$ on   
 $[0,T]\times L^2(\Fc_0;\R^d)$ which are continuously differentiable in $t$ $\in$ $[0,T)$, twice continuously Fr\'echet differentiable on  $L^2(\Fc_0;\R^d)$, and  which are liftings of functions on $[0,T]\times\Pc_{_2}(\R^d)$, i.e. 
 $\Phi(t,\xi)$ $=$ $\varphi(t,\P_\xi)$, for some $\varphi$ $\in$  $\Cc_b^{1,2}([0,T]\times\Pc_{_2}(\R^d))$, called inverse-lifted function of $\Phi$.

 \begin{Definition} \label{defvisco}
 We say that a locally bounded  function $u$ $:$ $[0,T]\times\Pc_{_2}(\R^d)$ $\rightarrow$ $\R$ is a viscosity (sub, super) solution to \reff{HJB2} if the lifted function $U$ $:$ $[0,T]\times L^2(\Fc_0;\R^d)$ $\rightarrow$ $\R$ defined by
 \beqs
 U(t,\xi) &=& u(t,\P_{_\xi}), \;\;\;\;\;  (t,\xi)  \in [0,T]\times L^2(\Fc_0;\R^d),
 \enqs
  is a viscosity (sub, super) solution to the lifted Bellman equation \reff{HJBlift}, that is: 
  
\noindent  (i) $U^*(T,.)$ $\leq$ $\hat G$, and for any test function $\Phi$ $\in$ $\Cc_\ell^2([0,T]\times L^2(\Fc_0;\R^d))$  such that  $U^*-\Phi$ has a maximum at $(t_{_0},\xi_{_0})$ $\in$ $[0,T)\times L^2(\Fc_0;\R^d)$, one has
 \beqs
-  \Dt{\Phi}(t_{_0},\xi_{_0}) +  \Hc(t_{_0},\xi_{_0},D\Phi(t_{_0},\xi_{_0}),D^2\Phi(t_{_0},\xi_{_0})) & \leq & 0. 
 \enqs 
(ii) $U_*(T,.)$ $\geq$ $\hat G$, and for any test function $\Phi$ $\in$ $\Cc_\ell^2([0,T]\times L^2(\Fc_0;\R^d))$  such that  $U_*-\Phi$ has a minimum at $(t_{_0},\xi_{_0})$ $\in$ $[0,T)\times L^2(\Fc_0;\R^d)$, one has
 \beqs
-  \Dt{\Phi}(t_{_0},\xi_{_0}) +   \Hc(t_{_0},\xi_{_0},D\Phi(t_{_0},\xi_{_0}),D^2\Phi(t_{_0},\xi_{_0})) & \geq & 0. 
 \enqs
 \end{Definition}

% where $\Cc_\mu(\R^d;\R^d)$ (resp. $\Cc_\mu(\R^d;\R^{d\times d})$) is the set of functions from $\R^d$ into  $\R^d$ (resp. $\R^{d\times d}$) continuous
 %on the support of $\mu$.

%We now turn to comparison principle (hence uniqueness result) for the Bellman equation \reff{HJB} that we rewrite in the ``Hamiltonian" form:
%\beq \label{HJB2}
%\Dt{v} + H(t,\mu,\partial_\mu v(t,\mu),\partial_x\partial_\mu v(t,\mu)) &=& 0    \;\;\; \mbox{ on } \; [0,T)\times\Pc_{_2}(\R^d),
%\enq
% where $H$ is the function defined  by
 %\beqs
 %H(t,\mu,p,\Gamma) &=& \inf_{\tilde\alpha\in L(\R^d;A)} \Big[ < f(t,.,\mu,\tilde\alpha(.),\tilde\alpha\star\mu) + p(.).b(t,.,\mu,\tilde\alpha(.),\tilde\alpha\star\mu) \\
% & & \hspace{2cm}  + \; \frac{1}{2}{\rm tr}\big(\Gamma(.)\sigma\sigma\trans(t,.,\tilde\alpha(.),\tilde\alpha\star\mu) \big), \mu > \Big],
 %\enqs
 %for  $(t,\mu)$ $\in$ $[0,T)\times\Pc_{_2}(\R^d)$, $(p,\Gamma)$ $\in$ $\Cc_\mu(\R^d;\R^d)\times\Cc_\mu(\R^d;\R^{d\times d})$, where
 %$\Cc_\mu(\R^d;\R^d)$ (resp. $\Cc_\mu(\R^d;\R^{d\times d})$) is the set of functions from $\R^d$ into  $\R^d$ (resp. $\R^{d\times d}$) continuous
 %on the support of $\mu$.

\vspace{1mm}

The main goal of this section is to prove the viscosity characterization  of the value function $v$ in \reff{defv} to the Bellman equation \reff{HJB}, hence equivalently the viscosity characterization of the 
lifted value function $V$ $:$ $[0,T]\times L^2(\Fc_0;\R^d)$ defined by 
\beqs
V(t,\xi) & = &  v(t,\P_{_\xi}), \;\;\; \xi \in L^2(\Fc_0;\R^d),
\enqs
to the  lifted Bellman equation \reff{HJBlift}.  We shall strenghten condition {\bf (H1)} by assuming in addition that $b,\sigma$ are uniformly continuous in $t$, and bounded  in $(a,\lambda)$:

\vspace{2mm}

\hspace{-7mm} {\bf (H1')}
%we have
%\beqs
%\int_0^T |b(t,0,\delta_0,0,\delta_0)|^2  +  |\sigma(t,0,\delta_0,0,\delta_0)|^2 dt & < & \infty.
%\enqs
There exists some constant $C_{b,\sigma}$ $>$ $0$  s.t. for all $t,t'$ $\in$ $[0,T]$, $x,x'$ $\in$ $\R^d$, $a,a'$ $\in$ $A$,
 $\lambda,\lambda'$ $\in$ $\Pc_{_2}(\R^d\times A)$,
\beqs
& & |b(t,x,a,\lambda) - b(t',x',a',\lambda')| + |\sigma(t,x,a,\lambda) - \sigma(t',x',a',\lambda')| \\
& \leq & C_{b,\sigma} \big[ m_{b,\sigma}(|t-t'|) + |x-x'| + |a-a'|  + W_2(\lambda,\lambda') \big],
\enqs
for some modulus  $m_{b,\sigma}$ (i.e. $m_{b,\sigma}(\tau)$ $\rightarrow$ $0$ when $\tau$ goes to zero) and
\beqs
|b(t,0,a,\delta_{_{0,0)}})|  +  |\sigma(t,0,a,\delta_{_{(0,0)}})| & \leq & C_{b,\sigma}.
\enqs

\vspace{1mm}

We also strenghten condition {\bf (H2)} by making additional (uniform) continuity  assumptions on the running and terminal cost functions, and 
boundedness conditions in $(a,\lambda)$:

\vspace{2mm}

\hspace{-7mm} {\bf (H2')}  (i) $g$ is  continuous on  $\R^d\times\Pc_{_2}(\R^d)$ and  there exists  some constant $C_{g}$ $>$ $0$ s.t.
for all  $x$ $\in$ $\R^d$, $\mu$ $\in$ $\Pc_{_2}(\R^d)$,
\beqs
|g(x,\mu)| & \leq & C_{g} \big(1 + |x|^2  + \|\mu\|^2_{_2} \big).
\enqs

\noindent  (ii)  There exists  some constant $C_{f}$ $>$ $0$ s.t. for all $t$ $\in$ $[0,T]$, $x$ $\in$ $\R^d$, $a$ $\in$ $A$,
$\lambda$ $\in$ $\Pc_{_2}(\R^d\times A)$,
\beqs
|f(t,x,a,\lambda)| & \leq & C_{f} \big(1 + |x|^2  + \|\lambda\|^2_{_2}\big),
\enqs
and some modulus  $m_{f}$ (i.e. $m_{f}(\tau)$ $\rightarrow$ $0$ when $\tau$ goes to zero) s.t. for all $t,t'$ $\in$ $[0,T]$, $x,x'$ $\in$ $\R^d$, $a$ $\in$ $A$,   $\lambda,\lambda'$ $\in$ $\Pc_{_2}(\R^d \times A)$,
\beqs
|f(t,x,a,\lambda) - f(t',x',a,\lambda')| & \leq & m_f\big(|t-t'| +  |x-x'|  + W_2(\lambda,\lambda')\big).  
\enqs 

\vspace{1mm}

The boundedness condition in {\bf (H1')}-{\bf (H2')} of $b,\sigma,f$ w.r.t. $(a,\lambda)$ $\in$ $A\times\Pc_{_2}(\R^d\times A)$ is ty\-pically satisfied when $A$ is bounded. 
Under {\bf (H1')}, we get  by standard arguments 
\beqs
\sup_{\alpha\in\Ac} \E\big[ \sup_{t\leq s\leq T} |X_s^{t,\xi}|^2 \big] & <  & \infty, 
\enqs
for any $t$ $\in$ $[0,T]$, $\xi$ $\in$ $L^2(\Fc_t;\R^d)$, which shows under the quadratic growth condition of $g$ and $f$ in {\bf (H2')} (uniformly in $a$) 
that  $v$ and $V$ also  satisfy  a quadratic growth condition: there exists some positive constant $C$ s.t.
\begin{equation}\label{vquadra}
\left\{
\begin{array}{ccl}
| v(t,\mu) | & \leq & C\big(1 + \|\mu\|_{_2}^2 \big), \;\;\; (t,\mu) \in [0,T] \times \Pc_{_2}(\R^d), \\
|V(t,\xi) | & \leq & C\big( 1 + \E|\xi|^2 \big), \;\;\; (t,\xi) \in [0,T]\times L^2(\Fc_0;\R^d),  
\end{array}
\right.
\end{equation}
and are  thus in particular locally bounded.

\vspace{2mm}

We first state a  flow continuity property of the marginal distribution of the controlled state process. Indeed, from  standard estimates on the state process under {\bf (H1')},
one easily checks (see also Lemma 3.1 in \cite{buetal14}) that  there exists some positive constant $C$, such that for all $\alpha$ $\in$ $\Ac$,
$t,t'$ $\in$ $[0,T$, $t\leq s\leq T, t'\leq s'\leq T$, $\mu=\P_\xi,\mu'=\P_{\xi'}$ $\in$ $\Pc_{_2}(\R^d)$:
\beqs 
\E\big| X_s^{t,\xi} - X_{s'}^{t',\xi'}\big|^2   & \leq & C\big(1 + \E|\xi|^2 +  \E|\xi'|^2 \big)\big( |t-t'| + |s-s'| +  \E|\xi-\xi'|^2 \big),
\enqs
and so from  the definition of the $2$-Wasserstein distance 
\beq \label{contiflow}
W_2(\P_s^{t,\mu},\P_{s'}^{t',\mu'}) & \leq & C\big(1 + \|\mu\|_{_2} +  \|\mu'\|_{_2}\big)\big( |t-t'|^{1\over 2} + |s-s'|^{1\over 2} +  W_2(\mu,\mu') \Big).
\enq

The next result states the viscosity property of the value function to the Bellman equation as a consequence of the dynamic programming principle \reff{DPPdeter}. 

\begin{Proposition}
The value function $v$ is a viscosity solution to the Bellman equation \reff{HJB}.
\end{Proposition}
{\bf Proof.}  We first show the continuity of $t$ $\mapsto$ $v(t,.)$ and $V(t,.)$  at $t$ $=$ $T$. For any $(t,\mu=\P_\xi)$ $\in$ $[0,T)\times\Pc_{_2}(\R^d)$, $\alpha$ $\in$ $\Ac$, we have from \reff{contiflow}
\beq \label{Wconti}
W_2(\P_T^{t,\mu},\mu) & \leq & \Big( \E|X_T^{t,\xi}-\xi|^2 \Big)^{1\over 2} \; \leq \; C(1 + \|\mu\|_{_2})|T-t|^{1\over 2},
\enq
for some positive constant $C$ (independent of $t$, $\mu$, $\alpha$). This means that $\P_T^{t,\mu}$ converges to $\mu$ in $\Pc_{_2}(\R^d)$ when
$t$ $\nearrow$ $T$,  uniformly in $\alpha$ $\in$ $\Ac$.  Now, from the definition of $v$ in \reff{defv}, we have
\beq
|v(t,\mu) - \hat g(\mu)| & \leq &  \sup_{\alpha\in\Ac}  \int_t^T \big| \hat f(s,\P_s^{t,\mu},\tilde\alpha(s,.,\P_s^{t,\mu})) \big| ds +  \big| \hat g(\P_T^{t,\mu}) - \hat g(\mu) \big| \nonumber \\
& \leq & C(1 +  \|\mu\|^2_{_2})|T-t| \; + \;   \sup_{\alpha\in\Ac} \big| \hat g(\P_T^{t,\mu}) - \hat g(\mu) \big|, \label{vgT}
\enq
from the growth condition on $f$ in {\bf (H2')}.  By the continuity assumption on $g$ together with the growth condition on $g$ in {\bf (H2')}, which allows to use dominated convergence  theorem,  we deduce from \reff{Wconti} that 
$\hat g(\P_T^{t,\mu})$ converges to $\hat g(\mu)$ when $t$ $\nearrow$ $T$,  uniformly in $\alpha$ $\in$ $\Ac$. This proves by \reff{vgT}  that $v(t,\mu)$ converges to $\hat g(\mu)$ when
$t$ $\nearrow$ $T$, i.e.  $v^*(T,\mu)$ $=$ $v_*(T,\mu)$ $=$ $\hat g(\mu)$ $=$ $v(T,\mu)$, and equivalently that $V(T,\xi)$ converges to $\hat G(\xi)$  when $t$ $\nearrow$ $T$, i.e.  
$V^*(T,\xi)$ $=$ $V_*(T,\mu)$ $=$ $\hat G(\xi)$ $=$ $V(T,\xi)$, 

Let us now prove the viscosity subsolution property of $V$ on  $[0,T)\times  L^2(\Fc_0;\R^d)$.  Fix $(t_{_0},\xi_{_0})$ $\in$ $[0,T)\times L^2(\Fc_0;\R^d)$, and consider some test function $\Phi$ $\in$ $\Cc_\ell^2([0,T]\times L^2(\Fc_0;\R^d))$ 
such that $V^*-\Phi$ has a maximum at $(t_{_0},\xi_{_0})$, and w.l.o.g. $V^*(t_{_0},\xi_{_0})$ $=$ $\Phi(t_{_0},\xi_{_0})$, so that $V^*$ $\leq$ $\Phi$. 
By definition of $V^*(t_{_0},\xi_{_0})$, there exists a sequence $(t_n,\xi_n)_n$ in $[0,T)\times L^2(\Fc_0;\R^d)$ s.t. 
\beqs
(t_n,\xi_n) \;\; \longrightarrow  \;\; (t_{_0},\xi_{_0}), & &  V(t_n,\xi_n) \; \longrightarrow \; V^*(t_{_0},\xi_{_0}),
\enqs
as $n$ goes to infinity.  By continuity of $\Phi$, we have
\beqs
\gamma_n \; := \; (V-\Phi)(t_n,\xi_n)  & \longrightarrow &  (V^* - \Phi)(t_{_0},\xi_{_0}) \; = \; 0,
\enqs
and let $(h_n)$ be a strictly positive sequence s.t. $h_n$ $\rightarrow$ $0$ and $\gamma_n/h_n$ $\rightarrow$ $0$. Consider the inverse-lifted function of $\Phi$, namely $\varphi$ $:$ $[0,T]\times\Pc_{_2}(\R^d)$ $\rightarrow$ $\R$ defined by $\varphi(t,\mu)$ $=$ $\Phi(t,\xi)$ for $t$ $\in$ $[0,T]$ and  $\mu$ $=$ $\P_{_\xi}$ $\in$ $\Pc_{_2}(\R^d)$, and recall that $\varphi$ $\in$  $\Cc_b^{1,2}([0,T]\times\Pc_{_2}(\R^d))$.    
Let $\tilde\alpha$ be an arbitrary element in  $L(\R^d;A)$, and consider the time-independent feedback control $\alpha$ $\in$ $\Ac$ associated with $\tilde\alpha$.  From the DPP \reff{DPPdeter} applied to $v(t_n,\mu_n)$, with $\mu_n$ 
$=$ $\P_{_{\xi_n}}$,  we have
\beqs
v(t_n,\mu_n)  & \leq &  \int_{t_n}^{t_n+h_n}  \hat f(s,\P_s^{t_n,\mu_n},\tilde\alpha)  ds \; + \;  v(t_n+h_n,\P_{t_n+h_n}^{t_n,\mu_n}).
\enqs
Since $v(t,\mu)$ $=$ $V(t,\xi)$ $\leq$ $V^*(t,\xi)$ $\leq$ $\Phi(t,\xi)$ $=$ $\varphi(t,\mu)$ for all $(t,\mu=\P_{_\xi})$ $\in$ $[0,T]\times\Pc_{_2}(\R^d)$,  this implies
%on a neighborhood of $(t_{_0},\mu_{_0})$, this  implies for $n$ large enough:
\beqs
\frac{\gamma_n}{h_n}  & \leq & \frac{1}{h_n}  \int_{t_n}^{t_n+h_n}  \hat f(s,\P_s^{t_n,\mu_n},\tilde\alpha)  ds \; + \;  \frac{\varphi(t_n+h_n,\P_{t_n+h_n}^{t_n,\mu_n}) - \varphi(t_n,\mu_n)}{h_n}.
\enqs
Applying It\^o's formula \reff{Ito} (similarly as in the verification theorem \ref{theoverif}) to $\varphi(s,\P_s^{t_n,\mu_n})$ between $t_n$ and $t_n+h_n$, we get
\beqs
\frac{\gamma_n}{h_n}  & \leq & \frac{1}{h_n}  \int_{t_n}^{t_n+h_n} \big[ \hat f(s,\P_s^{t_n,\mu_n},\tilde\alpha) +   \Dt{\varphi}(s,\P_s^{t_n,\mu_n}) + < \Lc_s^{\tilde\alpha}\varphi(s,\P_s^{t_n,\mu_n}), \P_s^{t_n,\mu_n}>  \big] ds
\enqs
Recall that $W_2(\mu_n,\mu_{_0})$ $\leq$ $\big(\E|\xi_n-\xi_0|^2\big)^{1\over 2}$, where $\mu_{_0}$ $=$ $\P_{_{\xi_0}}$, which shows that $\mu_n$ $\rightarrow$ $\mu_0$ in $\Pc_{_2}(\R^d)$ as $n$ goes to infinity. 
By the continuity of $b,\sigma,f,\varphi$ on their respective domains, the flow continuity property \reff{contiflow},  we then obtain by sending $n$ to infinity in the above inequality:
\beqs
0 & \leq & \hat f(t_{_0},\mu_{_0},\tilde\alpha) + \Dt{\varphi}(t_{_0},\mu_{_0}) + <\Lc_{t_{_0}}^{\tilde\alpha}\varphi(t_{_0},\mu_{_0}),\mu_{_0}>, 
\enqs
Since $\tilde\alpha$ is arbitrary in $L(\R^d;A)$, this shows 
\beqs
-  \Dt{\varphi}(t_{_0},\mu_{_0}) +  H(t_{0},\mu_{_0},\partial_\mu \varphi(t_{_0},\mu_{_0}),\partial_x\partial_\mu \varphi(t_{0},\mu_{_0})) & \leq &  0, 
\enqs
and thus at the lifted level: 
\beqs
-  \Dt{\Phi}(t_{_0},\xi_{_0}) +  \Hc(t_{0},\xi_{_0},D \Phi(t_{_0},\xi_{_0}), D^2\Phi(t_{0},\xi_{_0})) & \leq &  0, 
\enqs
which is the required  viscosity subsolution property.

We proceed finally with the viscosity supersolution property.  Fix $(t_{_0},\xi_{_0})$ $\in$ $[0,T)\times L^2(\Fc_0;\R^d)$, and consider some test function $\Phi$ $\in$ $\Cc_\ell^2([0,T]\times L^2(\Fc_0;\R^d))$ 
such that $V_*-\Phi$ has a minimum at $(t_{_0},\xi_{_0})$, and w.l.o.g. $V_*(t_{_0},\xi_{_0})$ $=$ $\Phi(t_{_0},\xi_{_0})$, so that $V_*$ $\geq$ $\Phi$. 
Again, by definition of $V_*(t_{_0},\xi_{_0})$, there exists a sequence $(t_n,\xi_n)_n$ in $[0,T)\times L^2(\Fc;\R^d)$ s.t.
$(t_n,\xi_n)$ $\longrightarrow$  $(t_{_0},\xi_{_0})$, and $V(t_n,\xi_n)$ $\longrightarrow$ $V_*(t_{_0},\xi_{_0})$ as $n$ goes to infinity.  We set $\gamma_n$ $:=$ $(V-\Phi)(t_n,\xi_n)$, which converges to zero, and we
consider a strictly positive sequence $(h_n)$ converging to zero and s.t. $\gamma_n/h_n$ also converges to zero.  Consider the inverse-lifted function of 
$\Phi$, namely $\varphi$ $\in$  $\Cc_b^{1,2}([0,T]\times\Pc_{_2}(\R^d))$ 
defined by $\varphi(t,\mu)$ $=$ $\Phi(t,\xi)$ for $t$ $\in$ $[0,T]$ and  $\mu$ $=$ $\P_{_\xi}$ $\in$ $\Pc_{_2}(\R^d)$. 
From the DPP \reff{DPPdeter},  for each $n$, and denoting by $\mu_n$ $=$ $\P_{_{\xi_n}}$ $\in$ $\Pc_{_2}(\R^d)$,  
there exists $\alpha^n$ $\in$ $\Ac$ associated to a feedback control $\tilde\alpha^n$ $\in$  $Lip([0,T]\times\R^d\times\Pc_{_2}(\R^d);A)$ s.t.
\beqs
v(t_n,\mu_n) + h_n^2  & \geq &  \int_{t_n}^{t_n+h_n}  \hat f(s,\P_s^{t_n,\mu_n},\tilde\alpha^n)  ds \; + \;  v(t_n+h_n,\P_{t_n+h_n}^{t_n,\mu_n}).
\enqs
Since $v(t,\mu)$ $=$ $V(t,\xi)$ $\geq$ $V_*(t,\xi)$ $\geq$ $\Phi(t,\xi)$ $=$ $\varphi(t,\mu)$ for all $(t,\mu=\P_{_\xi})$ $\in$ $[0,T]\times\Pc_{_2}(\R^d)$, this implies
\beqs
\frac{\gamma_n}{h_n} + h_n  & \geq & \frac{1}{h_n}  \int_{t_n}^{t_n+h_n}  \hat f(s,\P_s^{t_n,\mu_n},\tilde\alpha^n)  ds \; + \;  \frac{\varphi(t_n+h_n,\P_{t_n+h_n}^{t_n,\mu_n}) - \varphi(t_n,\mu_n)}{h_n}.
\enqs
Applying It\^o's formula \reff{Ito}  to $\varphi(s,\P_s^{t_n,\mu_n})$, we then get
\beqs
\frac{\gamma_n}{h_n} + h_n  & \geq & \frac{1}{h_n}  \int_{t_n}^{t_n+h_n} \big[   \Dt{\varphi}(s,\P_s^{t_n,\mu_n})   \; + \;   \hat f(s,\P_s^{t_n,\mu_n},\tilde\alpha^n(s,.,\P_s^{t_n,\mu_n}))  \\
& & \hspace{2cm} + \; < \Lc_s^{\tilde\alpha^n(s,.,\P_s^{t_n,\mu_n})}\varphi(s,\P_s^{t_n,\mu_n}), \P_s^{t_n,\mu_n}>  \big] ds  \\
& \geq &  \frac{1}{h_n}  \int_{t_n}^{t_n+h_n}
\Big(  \Dt{\varphi}(s,\P_s^{t_n,\mu_n}) \\
& & \hspace{1cm} + \;  \inf_{\tilde\alpha\in L(\R^d;A)} \big[ \hat f(s,\P_s^{t_n,\mu_n},\tilde\alpha) + < \Lc_s^{\tilde\alpha}\varphi(s,\P_s^{t_n,\mu_n}),\P_s^{t_n,\mu_n}> \big] \Big) ds.
\enqs
By sending $n$ to infinity together with the continuity assumption in {\bf (H1')}-{\bf (H2')} of $b,\sigma,f,\varphi$, uniformly in $a$ $\in$ $A$, and  the flow continuity property \reff{contiflow},  we get
%(which is uniform in $\alpha$ $\in$ $\Ac$), we get
\beqs
-  \Dt{\varphi}(t_{_0},\mu_{_0}) +  H(t_{0},\mu_{_0},\partial_\mu \varphi(t_{_0},\mu_{_0}),\partial_x\partial_\mu \varphi(t_{0},\mu_{_0})) & \geq &  0,
\enqs
which gives  the required  viscosity supersolution property of $V_*$, 
%i.e. 
%\beqs
%-  \Dt{\Phi}(t_{_0},\xi_{_0}) +  \Hc(t_{0},\xi_{_0},D \Phi(t_{_0},\xi_{_0}), D^2\Phi(t_{0},\xi_{_0})) & \geq &  0, 
%\enqs
and ends the proof.
\ep

\vspace{3mm}

%We now turn to comparison principle (hence uniqueness result) for the Bellman equation \reff{HJB} that we rewrite in the ``Hamiltonian" form:
%\beq \label{HJB2}
%\Dt{v} + H(t,\mu,\partial_\mu v(t,\mu),\partial_x\partial_\mu v(t,\mu)) &=& 0    \;\;\; \mbox{ on } \; [0,T)\times\Pc_{_2}(\R^d),
%\enq
% where $H$ is the function defined  by
 %\beqs
 %H(t,\mu,p,\Gamma) &=& \inf_{\tilde\alpha\in L(\R^d;A)} \Big[ < f(t,.,\mu,\tilde\alpha(.),\tilde\alpha\star\mu) + p(.).b(t,.,\mu,\tilde\alpha(.),\tilde\alpha\star\mu) \\
% & & \hspace{2cm}  + \; \frac{1}{2}{\rm tr}\big(\Gamma(.)\sigma\sigma\trans(t,.,\tilde\alpha(.),\tilde\alpha\star\mu) \big), \mu > \Big],
 %\enqs
 %for  $(t,\mu)$ $\in$ $[0,T)\times\Pc_{_2}(\R^d)$, $(p,\Gamma)$ $\in$ $\Cc_\mu(\R^d;\R^d)\times\Cc_\mu(\R^d;\R^{d\times d})$, where
 %$\Cc_\mu(\R^d;\R^d)$ (resp. $\Cc_\mu(\R^d;\R^{d\times d})$) is the set of functions from $\R^d$ into  $\R^d$ (resp. $\R^{d\times d}$) continuous
 %on the support of $\mu$.

We finally turn to   comparison principle (hence uniqueness result) for the Bellman equation \reff{HJB} (or \reff{HJB2}), hence equivalently for the lifted Bellman equation \reff{HJBlift}, which shall follow from 
comparison results for second order Hamilton-Jacobi equations in separable Hilbert space stated in \cite{lio89b}, see also \cite{fabgozswi15}.  We shall assume that  the $\sigma$-algebra  $\Fc_0$  is  
coun\-tably generated upto null sets, which ensures  that  the  Hilbert space $L^2(\Fc_0;\R^d)$ is separable, see \cite{doo94}, p. 92.  This is satisfied for example when $\Fc_0$ is the Borel $\sigma$-algebra 
of a canonical space $\Omega_0$ of continuous functions on $\R_+$,  in which case, $\Fc_0$ $=$ $\vee_{_{s\geq 0}} \Fc_s^{B^0}$, where $(\Fc_s^{B^0})$  is the canonical filtration on $\Omega_0$, and it is then known that 
$\Fc_0$ is countably generated, see for instance Exercise 4.21 in Chapter 1 of \cite{revyor99}.

%\marginpar{CHECK IF IT IS COHERENT WITH THE FACT THAT $\Fc_0$ is RICH ENOUGH} 

\begin{Proposition} \label{viscov}
Let $u$ and $w$ be two functions defined on  $[0,T]\times\Pc_{_2}(\R^d)$ satisfying a  quadratic growth  condition such that $u$ (resp. $w$)  is an upper (resp. lower) semicontinuous 
viscosity subsolution (resp. supersolution) to  \reff{HJB}. Then $u$ $\leq$ $w$. 
Consequently, the value function $v$ is the unique viscosity solution to the Bellman equation \reff{HJB} satisfying a quadratic growth condition \reff{vquadra}. 
\end{Proposition}
{\bf Proof.}  In view of our definition \ref{defvisco} of viscosity solution, we have to show a comparison principle for viscosity solutions to the lifted Bellman equation \reff{HJBlift}.   
We use the comparison principle proved in Theorem 3.50 in \cite{fabgozswi15} and only need to check that the hypotheses of this theorem are satisfied in our context for the lifted Hamiltonian $\Hc$ defined in \reff{defHc}.  Notice that the lifted Bellman equation \reff{HJBlift} is a bounded equation in the terminology 
of  \cite{fabgozswi15} (see their section  3.3.1) meaning that there is no linear dissipative operator on $L^2(\Fc_0;\R^d)$ in the equation. 
Therefore, the notion of $B$-continuity reduces to the standard notion of continuity in $L^2(\Fc_0;\R^d)$ since one can take for $B$ the identity operator. 
Their Hypothesis 3.44 follows from the uniform continuity  of $b$, $\sigma$, and $f$ in {\bf (H1')}-{\bf (H2')}.  Hypothesis 3.45  is immediately satisfied since there is no discount factor in our equation, i.e. $\Hc$ does not depend on $V$ but only on its derivatives.  The monotonicity condition in $Q$ $\in$ $S(L^2(\Fc_0;\R^d))$ of $\Hc$ in Hypothesis 3.46  is clearly satisfied.  Hypothesis 3.47  holds directly  when dealing with  bounded equations.   
Hypothesis 3.48  is obtained from the Lipschitz condition of $b,\sigma$ in {\bf (H1')}, and the uniform continuity condition on $f$ in {\bf (H2')}, while Hypothesis 3.49  follows from the quadratic growth condition of $\sigma$ in {\bf (H1')}. One can then apply Theorem 3.50  in  \cite{fabgozswi15} and conclude that comparison principle holds for the Bellman equation \reff{HJBlift}, hence for the Bellman equation \reff{HJB}.
\ep

\section{The case of open-loop controls}

\setcounter{equation}{0} \setcounter{Assumption}{0}
\setcounter{Theorem}{0} \setcounter{Proposition}{0}
\setcounter{Corollary}{0} \setcounter{Lemma}{0}
\setcounter{Definition}{0} \setcounter{Remark}{0}

In this section, we discuss how one can consider more generally open-loop controls instead of (Lipschitz) closed-loop  controls as imposed in the previous sections. 
%hence allowing a priori in particular bang-bang controls, which is useful in the applications. 
We shall restrict our framework  to usual controlled  McKean-Vlasov SDE with coefficients that do not depend on the law of the control but only on the law of the state process, hence in the form
\beq \label{McKeanopen}
dX_s &=& b(s,X_s,\alpha_s,\P_{_{X_s}}) ds +  \sigma(s,X_s,\alpha_s,\P_{_{X_s}}) dB_s,
\enq
where $b$, $\sigma$ are measurable  functions from $[0,T]\times\R^d\times A\times\Pc_{_2}(\R^d)$ into $\R^d$, respectively $\R^{d\times n}$, satisfying a Lipschitz condition: for all 
$t$ $\in$ $[0,T]$, $x,x'$ $\in$ $\R^d$, $a$ $\in$ $A$, $\mu,\mu'$ $\in$ $\Pc_{_2}(\R^d)$, 
\beq 
& & |b(t,x,a,\mu) - b(t,x',a,\mu')| + |\sigma(t,x,a,\mu) - \sigma(t,x',a,\mu')| \nonumber \\
& \leq & C \big[ |x-x'| + W_2(\mu,\mu') \big],   \label{lipopen}
\enq
for some positive constant $C$.  We denote by $\Ac_o$ the set of $\F$-progressive processes $\alpha$ valued in $A$, assumed for simplicity here to be a compact space of $\R^m$, 
%and satisfying the square-integrability condition: $\E[\int_0^T |\alpha_t|^2 dt]$ $<$ $\infty$, 
and consider the McKean-Vlasov control problem with open-loop controls when there is no running cost:
\beqs
\Vc_0 & := & \inf_{\alpha\in\Ac_o} \E \big[  g(X_T,\P_{_{X_T}}) \big]. 
\enqs
Under \reff{lipopen}, and given $t$ $\in$ $[0,T]$, $\xi$ $\in$ $L^2(\Fc_t;\R^d)$, $\alpha$ $\in$ $\Ac_o$, there exists a unique (pathwise and in law)  solution $X_s^{t,\xi}$ $=$ $X_s^{t,\xi,\alpha}$, $t\leq s\leq T$, solution to \reff{McKeanopen} starting from $\xi$ at time 
$t$, satisfying
\beqs
\E \big[ \sup_{t\leq s\leq T} |X_s^{t,\xi}|^2 \big] & \leq & C\big( 1 + \E|\xi|^2), 
% + \E\big[ \int_0^T |\alpha_t|^2 dt \big]  \big). 
\enqs
for some positive constant $C$ independent of $\alpha$ $\in$ $\Ac_o$. 
As in \reff{defPmu}, one can then define the flow $\P_s^{t,\mu}$ $=$ $\P_s^{t,\mu,\alpha}$, $t\leq s\leq T$, $\mu$ $\in$ $\Pc_{_2}(\R^d)$, 
$\alpha$ $\in$ $\Ac_o$, of the law of $X_s^{t,\xi}$, for $\mu$ $=$ $\P_\xi$, and it satisfies the flow property \reff{flowP}.  
We then define the value function in the Wasserstein space
\beq \label{vodet}
v_o(t,\mu) & := & \inf_{\alpha\in\Ac_o}\hat g(\P_T^{t,\mu}), \;\; t \in [0,T], \; \mu  \in \Pc_{_2}(\R^d), 
\enq
so that $\Vc_0$ $=$ $v_o(0,\P_{_{X_0}})$. Since the set of open-loop controls is larger than the  set of feedback controls, it is clear that $v_o$ is smaller than $v$ the value function to the McKean-Vlasov control problem with feedback controls considered  in the previous sections. By similar  arguments as in Theorem \ref{theoDPP},  one can show the DPP for the value function with open-loop controls, namely:
\beqs
v_o(t,\mu) &=& \inf_{\alpha\in\Ac_o}  v_o(\theta,\P_\theta^{t,\mu}), 
\enqs
for all $0\leq t\leq \theta\leq T$, $\mu$ $=$ $\P_\xi$  $\in$ $\Pc_{_2}(\R^d)$.  It would be possible to consider a nonzero running cost function $f$,  but in this case, one could not reformulate the 
value function $v_o$ as a deterministic control problem as in \reff{vodet}, and instead one has to consider the pair $(X_t,\P_{_{X_t}})$ as state variable in order to state a dynamic programming principle. This will be investigated 
in detail in  \cite{baycospha16}. 
From It\^o's formula \reff{Ito}, the infinitesimal version of the above DPP leads to the dynamic programming Bellman equation: 
\begin{equation} \label{HJBopen}
\left\{
\begin{array}{rcc}
- \partial_t v_o(t,\mu) + H_o\big(t,\mu,\partial_\mu v_o(t,\mu),\partial_x\partial_\mu v_o(t,\mu) \big) &=& 0, \;\;\; \mbox{ on } \; [0,T)\times\Pc_{_2}(\R^d), \\
v_o(T,.) &=& \hat g,  \;\;\; \mbox{ on }  \Pc_{_2}(\R^d)
\end{array}
\right. 
\end{equation}
where $H_o$ is the function defined by 
\beqs
H_o(t,\mu,p,\Gamma) &:=& - \inf_{\alpha\in\Ac_o} \E \big[  p(\xi).b(t,\xi,\alpha_t,\mu) +  
\frac{1}{2}{\rm tr}\big(\Gamma(\xi) \sigma\sigma\trans(t,\xi,\alpha_t,\mu) \big) \big], 
\enqs
for $(t,\mu)$ $\in$ $[0,T]\times\Pc_{_2}(\R^d)$, $(p,\Gamma)$ $\in$ $L_\mu^2(\R^d)\times L_\mu^\infty(\S^d)$, and with $\P_\xi$ $=$ $\mu$. 
Similarly as in Propositions \ref{theoverif} and \ref{viscov}, one can show a verification theorem for $v_o$ and prove that $v_o$ is the unique viscosity solution to \reff{HJBopen}. 

For any $\tilde\alpha$ $\in$ $L(\R^d;A)$, it is clear that the control $\alpha$ defined by $\alpha_s$ $=$ $\tilde\alpha(\xi)$, $t\leq s\leq T$, lies in $\Ac_o$, 
so that 
\beqs
H_o(t,\mu,p,\Gamma)  & \geq & - \inf_{\tilde\alpha\in L(\R^d;A)} \E \big[   p(\xi).b(t,\xi,\tilde\alpha(\xi),\mu)  +   \frac{1}{2}{\rm tr}\big(\Gamma(\xi) \sigma\sigma\trans(t,\xi,\tilde\alpha(\xi),\mu) \big) \big] \\
&   =  &   H(t,\mu,p,\Gamma),
\enqs
with $H$ the Hamiltonian in \reff{defH} for the McKean-Vlasov control problem with feedback control.  This inequality $H_o$ $\geq$ $H$ combined with  comparison principle for the Bellman equation \reff{HJBopen} is consistent with the inequality $v$ $\geq$ $v_o$. If we could prove that $H_o$ is equal to $H$ (which is not trivial in general), then this would show that $v_o$ is equal to $v$, i.e. the value functions to the McKean-Vlasov control problems with open-loop and feedback controls coincide. 
 Actually,  we notice that the minimization over the infinite dimensional space $\Ac_o$ in the Hamiltonian $H_o$ can be reduced into a minimization over the finite dimensional space $A$, namely:
\beq 
H_o(t,\mu,p,\Gamma) &=&  \tilde H_o(t,\mu,p,\Gamma) \label{Hointer}  \\
&:=& - < \inf_{a\in A} \big[   p(.).b(t,.,a,\mu) +  \frac{1}{2}{\rm tr}\big(\Gamma(.) \sigma\sigma\trans(t,.,a,\mu) \big) \big], \mu>.   \nonumber
\enq
Indeed, it is clear that $H_o$ $\leq$ $\tilde H_o$. Conversely, by continuity 
of the coefficients $b$, $\sigma$  w.r.t. the argument $a$ lying the compact space $A$, and invoking a measurable selection theorem, one can find for any 
$(t,\mu,p,\Gamma)$ $\in$ $[0,T]\times\Pc_{_2}(\R^d)\times L_\mu^2(\R^d)\times L_\mu^\infty(\S^d)$, 
 some measurable function $x$ $\in$ $\R^d$ $\mapsto$ $\hat a(t,x,\mu,p(x),\Gamma(x))$ $=$ $\hat\alpha(x)$ s.t.  for all $x$ $\in$ $\R^d$, 
\beqs
& & \inf_{a\in A} \big[  p(x).b(t,x,a,\mu) +  \frac{1}{2}{\rm tr}\big(\Gamma(x) \sigma\sigma\trans(t,x,a,\mu) \big) \big] \\
&=&  p(x).b(t,x,\hat\alpha(x),\mu) +  \frac{1}{2}{\rm tr}\big(\Gamma(x) \sigma\sigma\trans(t,x,\hat\alpha(x),\mu) \big). 
\enqs
By integrating w.r.t. $\mu$ $=$ $\P_\xi$, we then get
\beqs
\tilde H_o(t,\mu,p,\Gamma) &=& - \E \big[   p(\xi).b(t,\xi,\hat\alpha(\xi),\mu) +  
\frac{1}{2}{\rm tr}\big(\Gamma(\xi) \sigma\sigma\trans(t,\xi,\hat\alpha(\xi),\mu) \big) \big]  \\
& \leq &   H_o(t,\mu,p,\Gamma),
\enqs 
which shows the equality \reff{Hointer}.   Suppose now that there exists some smooth solution $w$ on $[0,T]\times\Pc_{_2}(\R^d)$ to the equation:
\begin{equation*}
\left\{
\begin{array}{rcc}
- \partial_t w(t,\mu) + \tilde H_o\big(t,\mu,\partial_\mu w(t,\mu),\partial_x\partial_\mu w(t,\mu) \big) &=& 0, \;\;\; \mbox{ on } \; [0,T)\times\Pc_{_2}(\R^d), \\
w(T,.) &=& \hat g,  \;\;\; \mbox{ on } \;  \Pc_{_2}(\R^d),
\end{array}
\right. 
\end{equation*}
such that for all $(t,\mu)$ $\in$ $[0,T)\times\Pc_{_2}(\R^d)$, the element $x$ $\mapsto$ 
$\hat a(t,x,\mu,\partial_\mu w(t,\mu)(x),\partial_x\partial_\mu w(t,\mu)(x))$ achieving the infimum in the definition of 
$\tilde H_o\big(t,\mu,\partial_\mu w(t,\mu),\partial_x\partial_\mu w(t,\mu) \big)$, is Lipschitz, i.e. lies in $L(\R^d;A)$, then 
(recall also Remark \ref{remoptcontrol}) 
\beqs
 \tilde H_o\big(t,\mu,\partial_\mu w(t,\mu),\partial_x\partial_\mu w(t,\mu) \big) &=& 
 H\big(t,\mu,\partial_\mu w(t,\mu),\partial_x\partial_\mu w(t,\mu) \big),
\enqs
which shows with \reff{Hointer} that $w$ solves both the Bellman equations \reff{HJBopen} and \reff{HJB2}.  By comparison principle, we conclude  
that $w$ $=$ $v$ $=$ $v_o$, which means in this case that the value functions to the McKean-Vlasov control problems with open-loop and feedback controls coincide.  
Such condition was satisfied for example in the case of the mean-variance portfolio selection problem studied  in paragraph \ref{subMV}.


\vspace{3mm}



\small

\begin{thebibliography}{}


\bibitem{ahmdin01} NU. Ahmed and X. Ding, Controlled McKean-Vlasov equation, {\it Communications in Applied Analysis}, {\bf 5}, (2001), 183-206.


\bibitem{ambetal05}  L. Ambrosio,  N.  Gigli  and G. Savar\'e, {\it Gradient Flows in Metric Spaces and in the Space of Probability Measures}, 
Lectures in Mathematics, Birkh\"auser Verlag, Basel (2005). 



\bibitem{anddje10} D. Andersson  and B. Djehiche, A maximum principle for SDEs of mean-field type, 
{\it Applied Mathematics and Optimization}, {\bf 63}, (2010), 341-356.


\bibitem{baycospha16} E. Bayraktar, A.  Cosso  and H. Pham, Randomized dynamic programming principle and Feynman-Kac representation for optimal control of McKean-Vlasov dynamics,   arXiv:1606.08204,  to appear in {\it Transactions of the American Mathematical Society}. 
%https://hal-auf.archives-ouvertes.fr/hal-01337515


%\bibitem{benetal13} Bensoussan A., Frehse J. and P.  Yam (2013): \emph{Mean Field Games and Mean Field Type Control Theory},
%Springer Briefs in Mathematics.


%\bibitem{bardap83} Barbu V.  and G. Da Prato (1983):  {\it Hamilton-Jacobi equations in Hilbert spaces}, Research Notes in Mathematics, n° 86, Pitman, Boston, MA, 1983.


\bibitem{benetal15} A. Bensoussan, J. Frehse  and P. Yam, The Master equation in mean-field theory, 
{\it Journal de Math\'ematiques Pures et Appliqu\'ees}, {\bf 103}, (2015), 1441-1474. 


\bibitem{benetal15b} A. Bensoussan, J. Frehse  and P. Yam, On the interpretation of the Master equation", arXiv: 1503.07754, 
to appear in {\it Stochastic Processes and their Applications}. 


\bibitem{benetal16} A. Bensoussan, K.C. Sung,  P. Yam  and  S.P. Yung, Linear-quadratic mean field games, 
{\it Journal of Optimization Theory and Applications},  {\bf 169}, (2016),  496-529.  



\bibitem{bjomur08} 
T. Bj\"ork, M.  Khapko and A.  Murgoci, On time inconsistent stochastic control in continuous time, arXiv:1612.03650, to appear in {\it Finance and Stochastics}. 
%Bj\"ork T. and A.  Murgoci (2008):  ``A general theory of Markovian time inconsistent stochastic control problems", working paper. 


\bibitem{bucetal11} R. Buckdahn, B. Djehiche  and J. Li, A general maximum principle for SDEs of mean-field type, {\it Applied Mathematics and Optimization}, {\bf 64}, (2011), 197-216. 


\bibitem{buetal14}  R. Buckdahn,  J. Li,  S. Peng  and  C. Rainer, Mean-field stochastic differential equations and associated PDEs, 
 arXiv: 1407.1215, to appear in the {\it Annals of Probability}.



 \bibitem{car12} P. Cardaliaguet,  Notes on mean field games, Notes from P.L. Lions lectures at Coll\`ege de France (2013). 
% https://www.ceremade.dauphine.fr/cardalia/MFG100629.pdf


\bibitem{cardel14b} R. Carmona  and F. Delarue, The Master equation for large population equilibriums,   
D. Crisan et al. (eds.), Stochastic Analysis and Applications 2014,  Springer Proceedings in Mathe\-matics $\&$  Statistics 100. 



\bibitem{cardel14} R. Carmona  and F. Delarue, Forward-backward Stochastic Differential Equations and Controlled McKean Vlasov Dynamics,  
 {\it Annals of Probability}, {\bf 43}, (2015), 2647-2700. 




\bibitem{cardel13} R. Carmona, F. Delarue  and A.  Lachapelle, Control of McKean-Vlasov dyna\-mics versus mean field games, 
 {\it Mathematics and Financial Economics}, {\bf  7}, (2013), 131-166.


\bibitem{caretal14}  R. Carmona, J.P.  Fouque  and   L. Sun, Mean field games and systemic risk,  {\it Communications
in Mathematical Sciences}, {\bf 13}, (2015), 911-933.  


 \bibitem{chacridel15} J.F. Chassagneux, D.   Crisan and F. Delarue, A probabilistic approach to classical solutions of the master equation for large population equilibria,  arXiv: 1411.3009. 


\bibitem{doo94} J.L. Doob, {\it Measure Theory},  Graduate texts in Mathematics, {\bf 143},  Springer (1994). 



\bibitem{fabgozswi15} G. Fabbri, F. Gozzi  and A. Swiech, {\it Stochastic optimal control in infinite dimension: dynamic programming and HJB equations} (with Chapter 6 by M. Fuhrman and G. Tessitore), (2015). 
%http://people.math.gatech.edu/~swiech/book.version11-10-2015.pdf



\bibitem{fenkat09}  J. Feng and M. Katsoulakis, A comparison principle for Hamilton-Jacobi equations related to controlled gradient flows in infinite dimensions,  {\it Archive for  Rational Mechanics and Analysis}, {\bf 192}, (2009),  275-310.


\bibitem{fisliv15} M. Fisher  and G. Livieri, Continuous time mean-variance portfolio optimization through the mean-field approach,  
{\it ESAIM Probability and Statistics}, {\bf 20}, (2016), 30-44. 




 \bibitem{FleSon06} W.H. Fleming  and H. M. Soner, \textit{Controlled Markov Processes and Viscosity Solutions}, 2nd edition,  Springer-Verlag (2006). 


\bibitem{ganetal08} W. Gangbo,  T.  Nguyen   and A. Tudorascu, Hamilton-Jacobi equations in the Wasserstein space, {\it Methods and Applications of Analysis},  {\bf 15}, (2008), 155-184.


\bibitem{ganswi15} W. Gangbo and A. Swiech, Metric viscosity solutions of Hamilton-Jacobi equations depending on local slopes,   
{\it Calculus of Variations and Partial Differential Equations}, {\bf 54}, (2015),  1183-1218.






\bibitem{huaetal06} M. Huang, P. Caines  and R. Malham\'e, Large population stochastic dynamic games: closed-loop McKean-Vlasov systems and the Nash certainty equivalence  principle, {\it Communication in Information and Systems}, {\bf 6}, (2006), 221-252.



\bibitem{jouetal08} B. Jourdain, S.  M\'el\'eard  and W. Woyczynski, Nonlinear SDEs driven by L\'evy processes and related PDEs,
{\it ALEA, Latin American Journal of Probability}, {\bf 4}, (2008),  1-29.


\bibitem{kac56} M. Kac, Foundations of kinetic theory, in Proceedings of the 3rd Berkeley Symposium on Mathematical Statistics and Probability, 
{\bf 3}, (1956), 171-197. 


\bibitem{laslio07} J.M. Lasry and  P.L. Lions, Mean-field games, {\it Japanese Journal of Mathematics}, {\bf 2}, (2007), 229-260. 


 \bibitem{laupir14} M. Lauri\`ere and O. Pironneau, Dynamic programming for mean-field type control, 
 %CRAS, {\bf 352}(9), 707-713.
{\it Journal of Optimization Theory and Applications}, {\bf 169}, (2016), 902-924. 

\bibitem{lizho00} D. Li  and X.Y. Zhou, Continuous-time mean-variance portfolio selection: a stochastic LQ framework,  
{\it Applied Mathematics and Optimization}, {\bf 42}, (2000), 19-33. 

\bibitem{lio88} P.L. Lions, Viscosity solutions of fully nonlinear second-order equations and optimal control in infinite dimension. Part I: the case of bounded stochastic evolution, 
{\it Acta Mathe\-matica}, {\bf 161}, (1988), 243-278. 

%\bibitem{lio89a} Lions P.L. (1989): ``Viscosity solutions of fully nonlinear second-order equations and optimal control in infinite dimension. Part II: Optimal control of Zakai's equation", 
%{\it Stochastic partial differential equations and applications, Lect. Notes in Math}, {\bf 1390}, Springer. 

\bibitem{lio89b} P.L. Lions, Viscosity solutions of fully nonlinear second-order equations and optimal control in infinite dimension. Part III: Uniqueness of viscosity solutions for general second-order equations,  {\it Journal of Functional Analysis}, {\bf 86}, (1989), 1-18. 





\bibitem{lio12} P.L. Lions, {\it Cours au Coll\`ege de France: Th\'eorie des jeux \`a champ moyens},  audio conference 2006-2012.  

 

%\bibitem{kur14}  Kurtz T. (2014):  ``Weak and strong solutions of general stochastic models", {\it  Electronic Communications in Probability}, {\bf 19}(58):1-16.


\bibitem{mckean67} H.P. McKean, Propagation of chaos for a class of nonlinear parabolic equations, Lecture Series in Differential Equations, {\bf 7}, 
(1967), 41-57. 


\bibitem{Pha09}  H. Pham,  \emph{Continuous-time stochastic control and applications with financial applications},  Series Stochastic Modeling and Applied Probability, {\bf 61}, Springer (2009). 


%\bibitem{phawei15} Pham H. and X. Wei (2015): ``Discrete time McKean-Vlasov control problem: a dynamic programming approach",  
%http://arxiv.org/pdf/1511.09273v1.pdf


\bibitem{revyor99} D. Revuz and M. Yor, {\it Continuous Martingales and Brownian Motion}, 3rd edition. New York, Berlin: Springer (1999). 


\bibitem{sni89} A.S. Sznitman, Topics in propagation of chaos, in {\it Lecture Notes in Mathematics}, Springer, {\bf 1464},  (1989), 165-251.


\bibitem{vil09} C. Villani, {\it Optimal transport, old and new}, Springer  (2009). 

\bibitem{yon13}  J. Yong, A linear-quadratic optimal control problem for mean-field stochastic diffe\-rential equations, {\it SIAM Journal on Control and Optimization}, {\bf 51}, (2013),  2809-2838.




\end{thebibliography}
\end{document}